\documentclass[12pt]{article}
\usepackage{amsmath,amsthm}
\usepackage{amssymb,latexsym}

\usepackage{enumerate}
\usepackage[english]{babel}
\usepackage{graphicx}

\usepackage{pdflscape}
\usepackage{array}
\usepackage[left=2.00cm, right=2.00cm, top=2.50cm, bottom=2.50cm]{geometry}

\def\Q{{\mathbb Q}}
\def\Z{{\mathbb Z}}

\newtheorem{lemma}{Lemma}
\newtheorem{theorem}[lemma]{Theorem}
\newtheorem{corollary}[lemma]{Corollary}

\title{
Integral bases and monogenity of composite fields
}
\author{
Istv\'{a}n Ga\'{a}l\thanks{
        Research supported in part by K115479 from the
        Hungarian National Foundation for Scientific Research
                         },\; 
and L\'aszl\'o Remete
\thanks{
        Research supported in part through the new 
				national excellence program of the ministry of human capacities.
                         }
\\ \\
University of Debrecen, Mathematical Institute \\
H--4002 Debrecen Pf.400., Hungary \\
e--mail: gaal.istvan@unideb.hu, remete.laszlo@science.unideb.hu
}

\begin{document}

\maketitle
\thispagestyle{empty}

\renewcommand{\thefootnote}{}



\renewcommand{\thefootnote}{\arabic{footnote}}
\setcounter{footnote}{0}

\begin{abstract}
We consider infinite parametric families of
high degree number fields composed of quadratic fields with pure 
cubic, pure quartic, pure sextic fields 
and with the so called simplest cubic, simplest quartic fields. 
We explicitly describe an integral basis of the composite fields.
We construct the index form, describe their factors and prove 
that the monogenity of the composite fields imply certain divisibility 
conditions on the parameters involved. These conditions usually can not hold, 
which implies the non-monogenity of the fields.

The fields that we consider are higher degree number fields, of degrees 6 up to 12.
The non-monogenity of the number fields is stated very often
as a consequence of the non-existence of the solutions of the index form equation.
Up to our knowlegde it is not at all feasible
to solve the index form equation in these high degree fields, especially not in
a parametric form.

On the other hand our method 
implies directly the non-monogenity in almost all cases.
We obtain our results in a parametric form, characterizing 
these infinite parametric families of composite fields.
\end{abstract}

\section{Introduction}

Monogenity of number fields and the existence of {\bf power integral 
bases} of type $\{1,\alpha,\ldots,\alpha^{n-1}\}$ is a classical
topic of algebraic number theory. The coefficients of the generators
of power integral bases are obtained as solutions of the
corresponding {\bf index form equations} cf. Section \ref{iiii}.

There are algorithms for the resolution of index form equations in
{\bf given specific} low degree number fields (for degrees 3 and 4 and some tedious
methods for degrees 5 and 6) and in some special type
of higher degree number fields see \cite{gaal}. 

We also succeeded to solve
the index form equation in certain {\bf infinite parametric families} of number fields,
using the solutions of a corresponding family of Thue equations or using
congruence considerations. 
Remark that very often we considered parametric families of number fields,
whose integral bases were not known in a parametric form. In those cases
we considered the problem of monogenity in the corresponding 
equation order.

We considered {\bf composites of number fields} of coprime discriminants
in \cite{comp1} showing a suitable factorization of the index form
in this case. In \cite{comp2}, \cite{comp3}, \cite{comp4} the authors
considered composites of  equation orders of number fields
and proved that under certain
congruence conditions on the defining polynomial these orders are not 
monogenic.

In some recent papers \cite{gr2}, \cite{gr4}
the authors developed a new and efficient technics to consider
monogenity in {\bf infinite parametric families of higher degree number fields}. 
The most important features of this method are that

\vspace{0.5cm}
\noindent
--the integral bases are determined in a parametric form, \\
--the factors of the index form are explicitly calculated,\\
--some linear combinations
of these factors are shown to have some non-trivial divisors.

\vspace{0.5cm}
In case a power integral basis (that is a solution of the
index form equation) exists, these imply some divisibility conditions on
the parameters. These divisibility conditions are usually not satisfied,
whence the fields are not monogenic.

Note that our algorithm to determine an integral basis
is based on standard methods cf. I.Stewart and D.Tall \cite{stewart}, J.Cook \cite{cook}.
We use this algorithm for parametric families of number fields. 
The explicit calculation of the factors of the index form
requires a very careful procedure because of the high degree, the
parameters and the large number of variables of the index form.

In our former results  \cite{gr2}, \cite{gr4}
we used this method to pure fields (up to degree 8)
and to the family of simplest sextic number fields.
Here we considerably extend our method and apply to 
composits of number fields up to degree 12. 
Up to our knowlegde this is the first time that
monogenity of such high degree fields are completely
characterized.

\vspace{0.5cm}

In the present paper {\bf we present an integral basis and obtain
conditions on the monogenity in composites of}
\begin{itemize}
\item quadratic fields and the simplest cubic fields (degree 6)
\item quadratic fields and pure cubic fields (degree 6)
\item quadratic fields and pure quartic fields (degree 8)
\item quadratic fields and the simplest quartic fields (degree 8)
\item the field $\Q(i\sqrt{3})$ and pure sextic fields (degree 12)
\end{itemize}
In each case we consider monogenity in the ring of integers of the composite 
field.

\section{Power integral bases and monogenity of number fields}
\label{iiii}

Here we shortly recall the concepts connected with monogenity of fields
\cite{gaal}, what we shall use throughout.
Let $\alpha$ be a primitive integral element of the number field $K$
(that is $K=\Q(\alpha)$) of degree $n$ with ring of integers $\Z_K$. 
The {\bf index} of $\alpha$ is
\[
I(\alpha)=(\Z_K^{+}:\Z[\alpha]^{+})=\sqrt{\left|\frac{D(\alpha)}{D_K}\right|}
=\frac{1}{\sqrt{|D_K|}}\prod_{1\leq i<j\leq n}\left|\alpha^{(i)}-\alpha^{(j)}\right| \;\; ,
\]
where $D_K$ is the discriminant of $K$ and $\alpha^{(i)}$ denote
the conjugates of $\alpha$. 

If $B=\{b_1=1,b_2,\ldots,b_n\}$ is an
integral basis of $K$, then the {\bf index form}
corresponding to this integral basis is defined by
\[
I(X_2,\ldots,X_n)=\frac{1}{\sqrt{|D_K|}}
\prod_{1\leq i<j\leq n}
\left(
(b_2^{(i)}-b_2^{(j)})X_2+\ldots+(b_n^{(i)}-b_n^{(j)})X_n
\right)
\]
(where $b_j^{(i)}$ denote the conjugates of $b_j$).
This is a homogeneous polynomial with integral coefficients
of degree $n(n-1)/2$. 
For the integral element 
\[
\alpha=x_1+b_2x_2+\ldots+b_nx_n
\]
(with $x_1,\ldots,x_n\in\Z$) we have
\[
I(\alpha)=|I(x_2,\ldots,x_n)|
\]
independently of $x_1$.
The element 
$\alpha$ generates a {\bf power integral basis} $\{1,\alpha,\ldots,\alpha^{n-1}\}$
if and only if $I(\alpha)=1$ that is $(x_2,\ldots,x_n)\in \Z^{n-1}$ is 
a solution of the
{\bf index form equation}
\begin{equation}
I(x_2,\ldots,x_n)=\pm 1 \;\;\; {\rm in} \;\;\; (x_2,\ldots,x_n)\in \Z^{n-1}.
\label{iiixxx}
\end{equation}
In this case the ring of integers of $K$ is a simple ring extension
of $\Z$, that is $\Z_K=\Z[\alpha]$ 
and $K$ is called {\bf monogenic}.

\vspace{1cm}
In our following statements and tables
for brevity we do not display the discriminants of the number fields
involved but they can be easily calculated from the discriminants
of the generating elements and the structure of the integral basis.
\vspace{1cm}

\section{Composites of \\quadratic fields and the simplest cubic fields}

Throughout this section we assume that
\begin{equation}
n,m \;\; {\rm are}\;\; {\rm integers}, n\not = 0,1, \;\;{\rm such}\;\;{\rm that} \;\;
n,m^2+3m+9 \;\;{\rm are}\;\;{\rm squarefree} \;\;{\rm and}\;\;{\rm coprime}.
\label{f1}
\end{equation}
In this section we consider the composite field $K=\Q(\alpha,\beta)$, where
\begin{center}
\begin{tabular}{rcl}
$\alpha$ &is a root of & $f(x)=x^2-n$,\\
$\beta$  &is a root of & $g(x)=x^3-mx^2-(m+3)x-1$.\\
\end{tabular}
\end{center}
The fields $M=\Q(\beta)$ are called {\em simplest cubic fields}, 
see \cite{shanks}.\\

\begin{theorem}
An integral basis of $K$ is given by
\[
\begin{array}{rl}
\displaystyle\left\{ 1, \beta, \beta^2, \frac{\alpha+1}{ 2}, \frac{\alpha\beta+\beta}{ 2}, \frac{\alpha\beta^2+\beta^2}{ 2} \right\}, & 
{\rm if}\;\; n\equiv 1\;(\bmod{\; 4}),\\ \\
\displaystyle\left\{1,\beta,\beta^2,\alpha,\alpha\beta,\alpha\beta^2\right\}, & 
{\rm if}\;\; n\equiv 2,3\;(\bmod{\; 4}).
\end{array}
\]
\end{theorem}

\vspace{0.5cm}

\noindent
{\bf Proof.}\\
It is well known that an integral basis $\{1,\omega\}$ and the discriminant  $D_L$ 
of $L=\Q(\alpha)$ are
\begin{equation}
\begin{array}{rl}
\omega=(1+\sqrt{n})/2,\; D_L=n
 & 
{\rm if}\;\; n\equiv 1\;(\bmod{\; 4}),\\ 
\omega=\sqrt{n},\; D_L=4n
 & 
{\rm if}\;\; n\equiv 2,3\;(\bmod{\; 4}).
\end{array}
\label{qu}
\end{equation}
More over, if $m^2+3m+9$ is squarefree then an  
integral basis $\{1,\delta_1,\delta_2\}$ and the discriminant  $D_M$ 
of $M=\Q(\beta)$ are
\[
\delta_1=\beta,\; \delta_2=\beta^2, \; D_M=(m^2+3m+9)^2.
\]
The discriminants $D_L,D_M$ are coprime, hence the composite field $K=LM$ has 
integral basis $\{1,\delta_1,\delta_2,\omega,\delta_1\omega,\delta_2\omega\}$ with
discriminant $D_K=D_L^3\cdot D_M^2$.
\mbox{}\\{$\Box$}

\vspace{0.5cm}

\begin{theorem}
If $K$ is monogenic then 
\[
\begin{array}{rl}
n\mid (m^2+3m+9) \pm1 \;\; {\rm and} \;\; (m^2+3m+9)\mid n^3\pm1,
 & 
{\rm if}\;\; n\equiv 1\;(\bmod{\; 4}),\\ \\
n\mid (m^2+3m+9) \pm1 \;\; {\rm and} \;\;   (m^2+3m+9)\mid 64n^3\pm1,
 & 
{\rm if}\;\; n\equiv 2,3\;(\bmod{\; 4}).
\end{array}
\]
\label{t2}
\end{theorem}

\vspace{0.5cm}

Throughout the paper the $\pm$ signs in the divisibility relations 
mean that the condition must hold either with $+$ or with $-$.

\vspace{0.5cm}

\noindent
{\bf Proof.}\\
The conjugates of $\alpha$ and $\beta$ are
\[
\alpha^{(1)}=\sqrt{n},\; \alpha^{(2)}=-\sqrt{n},\;
\beta^{(1)}=\beta,\;\beta^{(2)}={\displaystyle\frac{-1}{1+\beta}},\;
\beta^{(3)}={\displaystyle\frac{-1-\beta}{\beta}}.
\]

Set
\[
L^{(i,j)}=L^{(i,j)}(X_1,\ldots,X_6)=X_1+X_2\delta_1^{(j)}
+X_3\delta_2^{(j)}
+X_4\omega^{(i)}+X_5\omega^{(i)}\delta_1^{(j)}
+X_6\omega^{(i)}\delta_2^{(j)}.
\]
for $i=1,2;j=1,2,3$. Let
\begin{eqnarray}
F_1&=&\left( L^{(1,1)}-L^{(1,2)}\right)\left( L^{(1,1)}-L^{(1,3)}\right)\left( L^{(1,2)}-L^{(1,3)}\right)
\cdot \nonumber\\    
&&    \left( L^{(2,1)}-L^{(2,2)}\right)\left( L^{(2,1)}-L^{(2,3)}\right)\left( L^{(2,2)}-L^{(2,3)}\right)\nonumber \\
F_2&=&\left( L^{(1,1)}-L^{(2,1)}\right)\left( L^{(1,2)}-L^{(2,2)}\right)\left( L^{(1,3)}-L^{(2,3)}\right)\label{faktorok}\\
F_3&=&\left( L^{(1,1)}-L^{(2,2)}\right)\left( L^{(1,1)}-L^{(2,3)}\right)\left( L^{(1,2)}-L^{(2,1)}\right)
\cdot\nonumber \\    
&&    \left( L^{(1,2)}-L^{(2,3)}\right)\left( L^{(1,3)}-L^{(2,1)})(L^{(1,3)}-L^{(2,2)}\right). \nonumber
\end{eqnarray}
We find that 
\[
F_i(X_2,\ldots,X_6)=f_i\cdot G_i(X_2,\ldots,X_6) \;\; (i=1,2,3)
\]
where 
\[
f_1=m^2+3m+9=\sqrt{|D_M^2|},\;\; f_2=\sqrt{|D_L^3|},\;\; f_3=1
\]
and $G_i=G_i(X_2,\ldots,X_6) \; (i=1,2,3)$ are primitive polynomials with integer
coefficients. By 
\[
f_1f_2f_3=\sqrt{|D_K|},
\]
the index form equation corresponding to 
the given integral basis of $K$ is just
\[
G_1(x_2,\ldots,x_6)\cdot G_2(x_2,\ldots,x_6)\cdot G_3(x_2,\ldots,x_6)=\pm 1 
\;\; {\rm in} \;\; x_2,\ldots,x_6\in \Z.
\]
If $K$ admits a power integral basis, then there exist 
$x_2,\ldots,x_6\in\Z$ satisfying this equation, that is 
\[
G_i(x_2,\ldots,x_6)=\pm 1 \;\; (i=1,2,3),
\]
or equivalently 
\[
F_i(x_2,\ldots,x_6)=\pm f_i \;\; (i=1,2,3).
\]
Direct calculation of the factors show that the polynomials $F_1+F_3$ and $F_2^2-F_3$
have integer coefficients and
\[
n|F_1(X_2,\ldots,X_6)+F_3(X_2,\ldots,X_6) \;\; {\rm and} \;\; (m^2+3m+9)|F_2^2(X_2,\ldots,X_6)-F_3(X_2,\ldots,X_6).
\]
This immediately gives
\[
n|f_1\pm f_3 \;\; {\rm and}\;\;  (m^2+3m+9)|f_2^2\pm f_3.
\]
in case of a solution, which implies the assertion of Theorem \ref{t2}.
\mbox{}\\{$\Box$}

\vspace{1cm}

\section{Composites of \\quadratic and pure cubic fields}

Throughout this section we assume that
\begin{eqnarray}
&& n,m \;\; {\rm are}\;\; {\rm integers} \;\; n,m\not = 0,1, \nonumber\\
&&n \;\;{\rm is}\;\;{\rm squarefree}, \; m \;{\rm is}\; {\rm cubefree} 
\;\;{\rm and}\;\; \gcd(n,m)\in\{1,2,3\}, \label{f2}\\
&&m=uv^2 \; {\rm with} \; {\rm squarefree},\; {\rm integers}\; u,v,\;{\rm with}\;
(u,v)=1,\; {\rm and}\; 2,3\not| v.\nonumber
\end{eqnarray}
In this section we consider the composite field $K=\Q(\alpha,\beta)$, where
\begin{center}
\begin{tabular}{rcl}
$\alpha$ &is a root of & $f(x)=x^2-n$,\\
$\beta$  &is a root of & $g(x)=x^3-m$.\\
\end{tabular}
\end{center}

\begin{theorem}
According to the behaviour of $m,n$ mod 36, 
an integral basis of $K$ is given by the following table.
\label{th3}
\end{theorem}

\begin{center}
{\hspace*{-1cm}\noindent\scriptsize
\renewcommand{\arraystretch}{2.5}
\begin{tabular}{|m{2.5cm}|m{2.5cm}|c|}
\hline n mod 36 & m mod 36 & integral basis \\ 
\hline 1,5,13,17,25,29 & 1,10,17,19,26,35 & 
$\displaystyle\left\{ 1, \beta, \frac{\beta^2+uv^2\beta+v^2}{ 3v}, \frac{\alpha+1}{ 2}, \frac{(\alpha+1)\beta}{ 2}, \frac{(\alpha+1)(\beta^2+uv^2\beta+v^2}{ 6v}\right\} $  \\ 
\hline 1,5,13,17,25,29 & 2,3,5,6,7,11,13,14,\newline15,21,22,23,25,29,\newline30,31,33,34 & 
$\displaystyle\left\{ 1, \beta, \frac{\beta^2}{v}, \frac{\alpha+1}{ 2}, \frac{(\alpha+1)\beta}{ 2}, \frac{(\alpha+1)\beta^2}{ 2v} \right\}$ \\ 
\hline 2,7,10,11,14,19,22,\newline23,26,31,34,35 & 1,17,19,35 & 
$\displaystyle\left\{ 1, \beta, \frac{\beta^2+uv^2\beta+v^2}{ 3v}, \alpha,\alpha\beta, \frac{\alpha(\beta^2+uv^2\beta+v^2)}{3v}\right\}$ \\  
\hline 2,7,10,11,14,19,22,\newline23,26,31,34,35 & 3,5,7,11,13,15,21,\newline23,25,29,31,33 & $\displaystyle\left\{1,\beta,\frac{\beta^2}{v},\alpha,\alpha\beta,\frac{\alpha\beta^2}{v}\right\} $ \\
\hline 2,10,14,22,26,34 & 2,6,14,22,30,34 & $\displaystyle\left\{1,\beta,\frac{\beta^2}{v},\alpha,\alpha\beta,\frac{\alpha\beta^2}{2v}\right\} $ \\
\hline 2,10,14,22,26,34 & 10,26 & $\displaystyle\left\{ 1, \beta, \frac{\beta^2+uv^2\beta+v^2}{ 3v}, \alpha,\alpha\beta, \frac{\alpha(\beta^2+uv^2\beta+4v^2)}{ 6v} \right\} $ \\
\hline 7,11,19,23,31,35 & 2,6,14,22,30,34 & $\displaystyle\left\{1,\beta,\frac{\beta^2}{v},\alpha,\alpha\beta,\frac{(\alpha+1)\beta^2}{2v}\right\} $ \\
\hline 7,11,19,23,31,35 & 10,26 & $\displaystyle\left\{ 1, \beta, \frac{\beta^2+uv^2\beta+v^2}{ 3v}, \alpha,\alpha\beta, \frac{(\alpha+1)(\beta^2+uv^2\beta+4v^2)}{ 6v} \right\} $ \\
\hline 3,6,15,30 & 1,17,19,35 & $\displaystyle\left\{ 1, \beta, \frac{\beta^2+uv^2\beta+v^2}{ 3v}, \alpha,\frac{\alpha(\beta-u)}{3}, \frac{\alpha(\beta^2+uv^2\beta+v^2)}{ 3v} \right\} $ \\
\hline 3,6,15,30 & 5,7,11,13,23,25,29,31 & $\displaystyle\left\{1,\beta,\frac{\beta^2}{v},\alpha,\alpha\beta,\frac{\alpha(\beta^2+uv^2\beta+v^2)}{3v}\right\} $ \\
\hline 3,6,15,30 & 3,15,21,33 & $\displaystyle\left\{1,\beta,\frac{\beta^2}{v},\alpha,\alpha\beta,\frac{\alpha\beta^2}{3v}\right\} $ \\
\hline 
\end{tabular} }
\end{center}

\begin{center}
{\hspace*{-1cm}\noindent\scriptsize
\renewcommand{\arraystretch}{2.5}
\begin{tabular}{|m{2.5cm}|m{2.5cm}|c|}
\hline n mod 36 & m mod 36 & integral basis \\ 
\hline 21,33 & 1,10,17,19,26,35 & $\displaystyle\left\{ 1, \beta, \frac{\beta^2+uv^2\beta+v^2}{ 3v}, \frac{\alpha+1}{ 2}, \frac{(\alpha+3)(\beta-u)}{ 6}, \frac{(\alpha+1)(\beta^2+uv^2\beta+v^2)}{ 6v}\right\} $ \\
\hline 21,33 & 2,5,7,11,13,14, \newline 22,23,25,29,31,34 & $\displaystyle\left\{ 1, \beta, \frac{\beta^2}{v}, \frac{\alpha+1}{ 2}, \frac{(\alpha+1)\beta}{ 2}, \frac{(\alpha+3)(\beta^2+uv^2\beta+v^2)}{ 6v} \right\}$ \\
\hline 21,33 & 3,6,15,21,30,33 & $\displaystyle\left\{ 1, \beta, \frac{\beta^2}{v}, \frac{\alpha+1}{ 2}, \frac{(\alpha+1)\beta}{ 2}, \frac{(\alpha+3)\beta^2}{ 6v} \right\}$\\
\hline 6,30 & 10 & $\displaystyle\left\{ 1, \beta, \frac{\beta^2+v^2\beta+v^2}{ 3v}, \alpha,\frac{\alpha(\beta+2)}{3}, \frac{\alpha(\beta^2+2v^2\beta+6v^2)}{ 6v} \right\} $ \\
\hline 6,30 & 2,14,22,34 & $\displaystyle\left\{1,\beta,\frac{\beta^2}{v},\alpha,\alpha\beta,\frac{\alpha(\beta^2+uv^2\beta+4v^2)}{6v}\right\} $ \\
\hline 6,30 & 6,30 & $\displaystyle\left\{1,\beta,\frac{\beta^2}{v},\alpha,\alpha\beta,\frac{\alpha\beta^2}{6v}\right\} $ \\
\hline 6,30 & 26 & $\displaystyle\left\{1,\beta, \frac{\beta^2+2v^2\beta+v^2}{3v}, \alpha, \frac{\alpha(\beta+1)}{3}, \frac{\alpha(\beta^2+2v^2\beta+4v^2)}{ 6v} \right\} $ \\
\hline 3,15 & 10 & $\displaystyle\left\{ 1, \beta, \frac{\beta^2+v^2\beta+v^2}{ 3v}, \alpha,\frac{\alpha(\beta+2)}{3}, \frac{\alpha(\beta^2+2v^2\beta+6)+\beta^2+4v^2\beta+4v^2}{ 6v} \right\} $ \\
\hline 3,15 & 2,14,22,34 & $\displaystyle\left\{1,\beta,\frac{\beta^2}{v},\alpha,\alpha\beta,\frac{\alpha(\beta^2+uv^2\beta+4v^2)+3\beta^2}{6v}\right\} $ \\
\hline 3,15 & 6,30 & $\displaystyle\left\{1,\beta,\frac{\beta^2}{v},\alpha,\alpha\beta,\frac{(\alpha+3)\beta^2}{6v}\right\} $\\
\hline 3,15 & 26 & $\displaystyle\left\{1,\beta, \frac{\beta^2+2v^2\beta+v^2}{ 3v}, \alpha, \frac{\alpha(\beta+1)}{3}, \frac{(\alpha+1)(\beta^2+2v^2\beta+4v^2)}{ 6v} \right\} $ \\
\hline 
\end{tabular} }
\end{center}

\noindent
{\bf Proof.}\\
An integral basis $\{1,\omega\}$ and discriminant $D_L$ of $L=\Q(\alpha)$
are given in (\ref{qu}).
An integral basis $(1,\delta_1,\delta_2)$ and discriminant $D_M$ of $M=\Q(\beta)$
are given (cf. \cite{cohen}, Theorem 6.4.13) by:
\begin{equation}
\begin{array}{rrl}
\left\{1,\beta,\displaystyle{\frac{\beta^2}{v}}\right\},\;& D_M=-27u^2v^2 & 
{\rm if}\;\; u^2\not\equiv v^2\;(\bmod{\; 9}),\\ \\
\left\{1,\beta,\displaystyle{\frac{v^2+uv^2\beta+\beta^2}{3v}}\right\},\; & D_M=-3u^2v^2 & 
{\rm if}\;\; u^2\equiv v^2\;(\bmod{\; 9}).
\end{array}
\label{cub}
\end{equation}
Denote by $D_{K/L}$ and $D_{K/M}$ the relative discriminants of $K$ over $L$
and $M$, respectively. We have (cf. \cite{nark})
\begin{equation}
D_K=N_{L/\Q}(D_{K/L})\cdot D_L^3, \;\;\;\; D_K=N_{M/\Q}(D_{K/M})\cdot D_M^2.
\label{reldiscr}
\end{equation}
We have $\gcd (m,n)=1,2,3$. Denote by $\nu_p(k)$ the exponent of the prime $p$
in an integer $k$.
Set $u_0=u/(2^{\nu_2(u)}3^{\nu_3(u)})$, $n_0=n/(2^{\nu_2(n)}3^{\nu_3(n)})$.
The above discriminant relations imply
\[
3n_0^3(u_0v)^4|D_K.
\]
On the other hand $\{1,\delta_1,\delta_2,\omega,\delta_1\omega,\delta_2\omega\}$
is a basis in $K$ (not necessarily integral basis), whence $D_K$ divides 
the discriminant of this basis: 
\[
D_K|D_L^3D_M^2.
\]
Therefore $D_K$ must be of the form $\pm 2^r3^sn_0^3u_0^4v^4$.
Following the algorithm of \cite{cook} in order to obtain an integral basis, 
we start with the initial basis 
$\{b_1,\ldots,b_6\}=\{1,\delta_1,\delta_2,\omega,\delta_1\omega,\delta_2\omega\}$
and we test if its discriminant 
can be reduced by a 2-factor or by a 3-factor by interchanging one of its
elements by a new element. 

For $p=2$ and $p=3$ we perform the following procedure separately.
Let
\begin{equation}
\mu=
\frac{\lambda_1b_1+\ldots+\lambda_6b_6}{p}.
\label{lam}
\end{equation}
We let $\lambda_i \;(i=1,\ldots, 6)$ run through $\{0,1,\ldots,p-1\}$ and calculate the
defining polynomial of $\mu$:
\[
F(x)=\prod_{i=1}^6 \left(x- 
\frac{
\lambda_1b_1^{(i)}+\ldots+\lambda_6b_6^{(i)}}{p}
\right).
\]
Here $b_j^{(i)}$ denote the conjugates of $b_j$.
For each $\lambda_i\; (i=1,\ldots, 6)$ this polynomial can be written as
\[
F(x)=x^6+\frac{e_5}{p}x^5+\ldots+\frac{e_1}{p^5}x+\frac{e_0}{p^6}
\]
with integers $e_5\ldots,e_1,e_0$ depending on $m$ and $n$.
We calculated for which values $m (\bmod{p^6})$ and $n (\bmod{p^6})$
does this polynomial have integer coefficients, that is 
for which values $m (\bmod{p^6})$ and $n (\bmod{p^6})$ we have
\begin{eqnarray*}
e_5 &\equiv& 0 \; (\bmod{p}),\\
&\ldots &\\
e_1 &\equiv& 0 \; (\bmod{p^5}),\\
e_0 &\equiv& 0 \; (\bmod{p^6}).
\end{eqnarray*}
For the selected pairs $m (\bmod{p^6})$ and $n (\bmod{p^6})$
the basis element having coefficient $\lambda_i=1$
can be replaced by the above element $\mu$ to diminish the discriminant by a $p^2$ factor.
(Any non-zero $\lambda_i$ can be transformed into $\lambda_i=1$ by multiplying
(\ref{lam}) by the inverse of $\lambda_i$ modulo $p$ and by subtracting a suitable 
integer element.)

This procedure is continued with $p$ until no reduction of the 
discriminant of the basis is possible. Then the same procedure is 
performed for the other value of $p$, as well.

Finally we combine the basis $\{b_1,\ldots,b_6\}$ of the 2-maximal order
and the basis $\{f_1,\ldots,f_6\}$ of the 3-maximal order of $K$
into a basis of the maximal order of $K$.

First we remark that we choose the basis elements $b_i$ and $f_i$
so that in the numerator of $b_i$ and $f_i$
the coefficient of
$x^{i-1}$ is equal to 1 
and the coefficients of $x^i,x^{i+1},\ldots$ are equal to 0.
We construct 
the integral basis $h_1=1,h_2,\ldots,h_n$ of $K$ with the same property.
Assume that $b_i$ has denominator $2^{k_2}$ and 
$f_i$ has denominator $3^{k_3}$.
Calculate $y_2$ and $y_3$ with 
$2^{k_2}y_2\equiv 1(\bmod{\; 3^{k_3}})$ and
$3^{k_3}y_3\equiv 1(\bmod{\; 2^{k_2}})$. Then $x=2^{k_2}y_2+3^{k_3}y_3$ is a solution of
the system $x\equiv 1 (\bmod{\; 2^{k_2}})$, $x\equiv 1 (\bmod{\; 3^{k_3}})$. We also have
$x\equiv 1 (\bmod{\; 2^{k_2}3^{k_3}})$. We set 
\[
h_i'=y_3b_i+y_2f_i=\frac{y_3 3^{k_3} (2^{k_2}b_i) + y_2 2^{k_2} (3^{k_3}f_i)   }{2^{k_2}3^{k_3}}
\]
which is an algebraic integer. The coefficient of $x^{i-1}$ in $h_i'$ is
\[
c_{i-1}=\frac{y_3 3^{k_3} + y_2 2^{k_2} }{2^{k_2}3^{k_3}}.
\]
The numerator can be written as $1+\ell\cdot  2^{k_2}3^{k_3}$ (with an integer $\ell$), hence 
\[
c_{i-1}=\frac{1}{2^{k_2}3^{k_3}}+\ell.
\]
We set $h_i=h_i'-\ell x^{i-1}$ which is also an algebraic integer.
In the numerator of $h_i$ the coefficient of $x^{i-1}$ is equal to 1 and
the coefficients of $x^i,x^{i+1},\ldots$ are equal to 0.

We show that $b_i$ can be expressed as a linear combination of $(1,x,\ldots,x^{i-2},h_i)$
(the same holds for $f_i$) which proves that $\{h_1=1,h_2,\ldots,h_n\}$ is indeed an 
integral basis, being a 2-maximal order and a 3-maximal order of $K$.

We make use of $3^{k_3}y_3 \equiv 1 (\bmod{\; 2^{k_2}})$ here, which implies 
$3^{k_3}y_3=1+q 2^{k_2}$ with an integer $q$.
We have 
\[
3^{k_3}h_i'=y_33^{k_3}b_i+y_23^{k_3}f_i=(1+q 2^{k_2})b_i+y_23^{k_3}f_i=
b_i+q 2^{k_2}b_i+y_23^{k_3}f_i.
\]
Observe that $2^{k_2}b_i$ and $3^{k_3}f_i$ is a linear combination of 
$(1,x,\ldots,x^{i-1})$ with integer coefficients.
We have
\[
b_i\in {\cal L}(1,x,\ldots,x^{i-1},h_i')\subseteq{\cal L}(1,x,\ldots,x^{i-2},h_i)
\]
since the coefficient of $x^{i-1}$ in the numerator of $h_i$ is 
equal to 1. 
($\cal L$ denotes the set of the linear combinations with integer coefficients
of the elements involved.)
This implies our assertion.

Finally we collect those pairs $m,n$
for which we obtained the same type of integral basis.
\mbox{}\\{$\Box$}

\vspace{0.5cm}

\noindent
{\bf Remark.} The above method obviously works for any distinct primes, as well.

\begin{theorem}
If $K$ admits a power integral basis, then the following divisibility 
conditions must hold:
\end{theorem}

\[
\begin{array}{|c|c|c|c|}
\hline
\begin{array}{c}\\n \; {\rm mod} \; 36 \\ \mbox{}  \end{array}  
& m \; {\rm mod} \; 36  & 1 &  2 \\ \hline
1,5,13,17,25,29 & 1,10,17,19,26,35& n\mid 3m^2 \pm1 & m\mid n^3\pm1 \\ \hline
1,5,13,17,25,29 & \begin{array}{l}2,3,5,6,7,11,13,14,\\ 15,21,22,23,25,29,\\ 30,31,33,34 \end{array}&
n\mid 27m^2 \pm1 & 9m\mid n^3\pm1\\ \hline
\begin{array}{l}2,7,10,11,14,19,22,\\23,26,31,34,35\end{array}&
1,17,19,35&4n\mid 3m^2 \pm1 & m\mid 64n^3\pm1\\ \hline
\begin{array}{l}2,7,10,11,14,19,22,\\23,26,31,34,35\end{array}&
\begin{array}{l}3,5,7,11,13,15,21,\\23,25,29,31,33\end{array}&
4n\mid 27m^2 \pm1 & 9m\mid 64n^3\pm1\\ \hline
\begin{array}{l}2,7,10,11,14,19,\\22,23,26,31,34,35\end{array}
 & 2,6,14,22,30,34 &
\displaystyle 4n\mid \frac{27m^2}{2} \pm2 & 9m\mid 16n^3\pm2\\ \hline
\begin{array}{l}2,7,10,11,14,19,\\22,23,26,31,34,35\end{array}&
10,26&\displaystyle 4n\mid \frac{3m^2}{2} \pm2 & m\mid 16n^3\pm2
\\ \hline
3,6,15,30&1,17,19,35&\displaystyle \frac{4n}{3}\mid m^2 \pm3 & \displaystyle 3m\mid \frac{64n^3}{9}\pm3
\\ \hline
3,6,15,30&\begin{array}{l}3,5,7,11,13,15,\\21,23,25,29,31,33\end{array}&
4n\mid 9m^2 \pm3 & \displaystyle 9m\mid \frac{64n^3}{9}\pm3
\\ \hline
21,33 & 1,10,19&
\displaystyle \frac{n}{3}\mid m^2 \pm3 & \displaystyle 3m\mid \frac{n^3}{9}\pm3
\\ \hline
21,33 &
\begin{array}{l}2,3,5,6,7,11,13,\\14,15,21,22,23,25,\\29,30,31,33,34\end{array}
&\displaystyle n\mid 9m^2 \pm3&\displaystyle 9m\mid \frac{n^3}{9}\pm3
\\ \hline
21,33 & 17,26,35 & 
\displaystyle \frac{n}{3}\mid m^2 \pm3 & \displaystyle m\mid \frac{n^3}{9}\pm3
\\ \hline
3,6,15,30&10,26&
\displaystyle \frac{4n}{3}\mid \frac{m^2}{2} \pm6 & \displaystyle 3m\mid \frac{16n^3}{9}\pm6
\\ \hline
3,6,15,30&2,6,14,22,30,34&
\displaystyle 4n\mid \frac{9m^2}{2} \pm6 & \displaystyle 9m\mid \frac{16n^3}{9}\pm6
\\ \hline
\end{array}
\]

\vspace{1cm} 

\noindent
{\bf Proof.}\\
Let 
\[
\alpha^{(1)}=\sqrt{n},\;\alpha^{(2)}=-\sqrt{n}\;
\beta^{(1)}=\sqrt[3]{m},\;\beta^{(2)}=\varepsilon\sqrt[3]{m},\;
\beta^{(3)}=\varepsilon^2\sqrt[3]{m}
\]
with $\varepsilon=\exp(2\pi i/3)$. 
Assume that $\{b_1=1,b_2,\ldots,b_6\}$ is an integral basis of $K$.
The elements of the integral basis are composed of 
$\alpha$ and $\beta$ hence it is unique to denote the conjugate of $b_k$ 
corresponding to $\alpha^{(i)}$ and $\beta^{(j)}$ by $b_k^{(i,j)}$ 
($i=1,2;j=1,2,3$). Set
\[
L^{(i,j)}=L^{(i,j)}(X_1,X_2,\ldots,X_6)=
X_1+b_2^{(i,j)}X_2+\ldots+b_6^{(i,j)}X_6.
\]
Using the $L^{(i,j)}$ we construct the same polynomials $F_1,F_2,F_3$
like in (\ref{faktorok}). We write these polynomials in the form
\[
F_i(X_2,\ldots,X_6)=f_i\cdot G_i(X_2,\ldots,X_6) \;\; (i=1,2,3)
\]
where $f_i$ are integers or square roots of integers with $(f_1f_2f_3)^2=|D_K|$
depending on the parameters only and $G_i$ are primitive polynomials with
integer coefficients.
Then the index form equation corresponding to the basis
$\{b_1=1,b_2,\ldots,b_6\}$ can be written as 
\[
G_1(x_2,\ldots,x_6)\cdot G_2(x_2,\ldots,x_6)\cdot G_3(x_2,\ldots,x_6)=\pm 1,
\;\; {\rm in}\;\; x_2,\ldots,x_6\in\Z,
\]
that is
\[
G_i(x_2,\ldots,x_6)=\pm 1 \;\; (i=1,2,3),
\]
or equivalently
\[
F_i(x_2,\ldots,x_6)=\pm f_i \;\; (i=1,2,3).
\]
If a power integral basis in $K$ exists then there exist $x_2,\ldots,x_6\in\Z$
satisfying the above equations.
Calculating the explicit form of the $F_i$ in each case of the integral basis
we always have 
\[
n|F_1(X_2,\ldots,X_6)+F_3(X_2,\ldots,X_6) \;\; {\rm and} \;\; m|F_2^2(X_2,\ldots,X_6)-F_3(X_2,\ldots,X_6),
\]
where again the polynomials $F_1+F_3$ and $F_2^2-F_3$
have integer coefficients.
These imply the divisibility conditions of the Theorem.
\mbox{}\\{$\Box$}

\vspace{1cm}

\section{Composites of \\quadratic and pure quartic fields}

Throughout this section we assume that
\begin{eqnarray}
&&n,m \;\; {\rm are}\;\; {\rm squarefree}\;\; {\rm integers}, n,m\not = 0,1, \nonumber\\
&& \gcd(n,m)\in\{1,2\},\label{f3}\\
&& {\rm if} \; m=n, \; {\rm then}\; m=n\not\in\{-2,-1,2\}. \nonumber
\end{eqnarray}
In this section we consider the composite field $K=\Q(\alpha,\beta)$, where
\begin{center}
\begin{tabular}{rcl}
$\alpha$ &is a root of & $f(x)=x^2-n$,\\
$\beta$  &is a root of & $g(x)=x^4-m$.\\
\end{tabular}
\end{center}

\begin{theorem}
According to the behaviour of $m,n$ mod 8, 
an integral basis of $K$ is given by the following table:
\label{hhhh}
\end{theorem}

{\noindent\tiny
\renewcommand{\arraystretch}{2.5}
\begin{tabular}{|c|c|c|}
\hline n mod 8 & m mod 8 & integral basis \\ 
\hline 1,5 & 2,3,6,7 & 
$\displaystyle\left\{ 1, \beta, \beta^2, \beta^3, \frac{\alpha+1}{ 2}, \frac{\beta\alpha+\beta}{ 2}, \frac{\beta^2\alpha+\beta^2}{ 2}, \frac{\beta^3\alpha+\beta^3}{ 2}\right\} $ \\ 
\hline 1,5 & 1 & 
$\displaystyle\left\{ 1, \beta, \frac{\beta^2+1}{ 2}, \frac{\beta^3+\beta^2+\beta+1}{ 4}, \frac{\alpha+1}{ 2}, \frac{\beta\alpha+\beta}{ 2}, \frac{\beta^2\alpha+\beta^2+\alpha+1}{ 4}, \frac{\beta^3\alpha+\beta^3+\beta^2\alpha+\beta^2+\beta\alpha+\beta+\alpha+1}{ 8} \right\} $\\ 
\hline 1,5 & 5 & 
$\displaystyle\left\{ 1, \beta, \frac{\beta^2+1}{ 2}, \frac{\beta^3+\beta}{ 2}, \frac{\alpha+1}{ 2}, \frac{\beta\alpha+\beta}{ 2}, \frac{\beta^2\alpha+\beta^2+\alpha+1}{ 4}, \frac{\beta^3\alpha+\beta^3+\beta\alpha+\beta}{ 4} \right\} $\\ 
\hline 2,6 & 1 & 
$\displaystyle\left\{
1, \beta, \frac{\beta^2+1}{ 2}, \frac{\beta^3+\beta^2+\beta+1}{ 4}, \alpha, \frac{\beta\alpha+\alpha}{ 2}, \frac{\beta^2\alpha+\alpha}{ 2}, \frac{\beta^3\alpha+\beta^2\alpha+\beta\alpha+\alpha}{ 4} \right\} $\\ 
\hline 2,6 & 5 & 
$\displaystyle\left\{1, \beta, \frac{\beta^2+1}{ 2}, \frac{\beta^3+\beta}{ 2}, \alpha, \frac{\beta\alpha+\alpha}{ 2}, \frac{\beta^2\alpha+\alpha}{ 2}, \frac{\beta^3\alpha+\beta^2\alpha+\beta\alpha+\alpha}{ 4}
\right\} $\\
\hline 3,7 & 2,6 & 
$\displaystyle\left\{1, \beta, \beta^2, \beta^3, \alpha, \frac{\beta^3+\beta\alpha+\beta}{ 2}, \frac{\beta^2\alpha+\beta^2}{ 2}, \frac{\beta^3\alpha+\beta^3}{ 2}
\right\} $\\
\hline 3,7 & 3,7 & 
$\displaystyle\left\{1, \beta, \beta^2, \beta^3, \frac{\beta^2+\alpha}{ 2}, \frac{\beta^3+\beta\alpha}{ 2}, \frac{\beta^2\alpha+1}{ 2}, \frac{\beta^3\alpha+\beta}{ 2}\right\} $\\
\hline 3,7 & 1 & 
$\displaystyle\left\{1, \beta, \frac{\beta^2+1}{ 2}, \frac{\beta^3+\beta^2+\beta+1}{ 4}, \alpha, \frac{\beta\alpha+\beta+\alpha+1}{ 2}, \frac{\beta^2\alpha+\beta^2+2\beta+3\alpha+1}{ 4}, \frac{\beta^3\alpha+\beta^2\alpha+\beta\alpha+\alpha}{ 4}\right\} $\\
\hline 3,7 & 5 & 
$\displaystyle\left\{1, \beta, \frac{\beta^2+1}{ 2}, \frac{\beta^3+\beta}{ 2}, \alpha, \frac{\beta\alpha+\beta+\alpha+1}{ 2}, \frac{\beta^2\alpha+\alpha}{ 2}, \frac{\beta^3\alpha+\beta^3+\beta^2\alpha+\beta^2+\beta\alpha+\beta+\alpha+1}{ 4}\right\} $\\
\hline 2 & 2 & 
$\displaystyle\left\{1, \beta, \beta^2, \beta^3, \frac{\beta^2+\alpha}{ 2}, \frac{\beta^3+\beta\alpha}{ 2}, \frac{\beta^2\alpha+2}{ 4}, \frac{\beta^3\alpha+2\beta}{ 4}\right\} $\\
\hline 2 & 3 & 
$\displaystyle\left\{1, \beta, \beta^2, \beta^3, \frac{\beta^2\alpha+\alpha}{ 2}, \frac{\beta^2+\beta\alpha+1}{ 2}, \frac{\beta^3+\beta^2\alpha+\beta}{ 2}, \frac{\beta^3\alpha+2\beta^3+\beta^2\alpha+\beta\alpha+\alpha+2}{ 4}\right\} $\\
\hline 2 & 6 & 
$\displaystyle\left\{1, \beta, \beta^2, \beta^3, \frac{\beta^2+\alpha}{ 2}, \frac{\beta^3+\beta\alpha}{ 2}, \frac{\alpha\beta^2}{ 2}, \frac{\beta^3\alpha+2\beta^3+2\beta}{ 4}\right\} $\\
\hline 2 & 7 & 
$\displaystyle\left\{1, \beta, \beta^2, \beta^3, \frac{\beta^2\alpha+2\beta+\alpha}{ 4}, \frac{\beta^2+\beta\alpha+1}{ 2}, \frac{\beta^3+\beta^2\alpha+\beta}{ 2}, \frac{\beta^3\alpha+2\beta^2+\beta\alpha}{ 4}\right\} $\\
\hline 6 & 2 & 
$\displaystyle\left\{1, \beta, \beta^2, \beta^3, \frac{\beta^2+\alpha}{ 2}, \frac{\beta^3+\beta\alpha}{ 2}, \frac{\alpha\beta^2}{ 2}, \frac{\beta^3\alpha+2\beta^3+2\beta}{ 4}\right\} $\\
\hline 6 & 3 & 
$\displaystyle\left\{1, \beta, \beta^2, \beta^3, \frac{\beta^2\alpha+\alpha}{ 2}, \frac{\beta^2+\beta\alpha+1}{ 2}, \frac{\beta^3+\beta^2\alpha+\beta}{ 2}, \frac{\beta^3\alpha+\beta^2\alpha+2\beta^2+\beta\alpha+2\beta+\alpha}{ 4}
\right\} $\\
\hline 6 & 6 & 
$\displaystyle\left\{
1, \beta, \beta^2, \beta^3, \frac{\beta^2+\alpha}{ 2}, \frac{\beta^3+\beta\alpha}{ 2}, \frac{\beta^2\alpha+2}{ 4}, \frac{\beta^3\alpha+2\beta}{ 4}\right\} $\\
\hline 6 & 7 & 
$\displaystyle\left\{1, \beta, \beta^2, \beta^3, \frac{2\beta^3+\beta^2\alpha+\alpha}{ 4}, \frac{\beta^2+\beta\alpha+1}{ 2}, \frac{\beta^3+\beta^2\alpha+\beta}{ 2}, \frac{\beta^3\alpha+\beta\alpha+2}{ 4}\right\} $\\
\hline 
\end{tabular} }\\

\vspace{0.5cm}

{\bf Proof.}\\
An integral basis $\{1,\omega\}$ and discriminant $D_L$ of $L=\Q(\alpha)$
are given in (\ref{qu}).
An integral basis $\{1,\delta_2,\delta_3,\delta_4\}$ 
and discriminant $D_M$ of $M=\Q(\beta)$
are given by (cf. \cite{gr2}):
\begin{equation*}
\begin{array}{rrl}
\{1,\beta,\beta^2,\beta^3\},\;& D_M=-256m^3,& {\rm if}\;\; m\equiv 2,3\;(\bmod{\; 4}),\\ \\
\left\{1,\beta,\displaystyle{\frac{1+\beta^2}{2}},
\displaystyle{\frac{1+\beta+\beta^2+\beta^3}{4}}\right\},
\; & D_M=-4m^3, & {\rm if}\;\; m\equiv 1\;(\bmod{\; 8}),
\\  \\
\left\{1,\beta,\displaystyle{\frac{1+\beta^2}{2}},
\displaystyle{\frac{\beta+\beta^3}{2}}\right\},
\; & D_M=-16m^3, & {\rm if}\;\; m\equiv 5\;(\bmod{\; 8}).
\end{array}
\end{equation*}
Starting from the initial basis 
$\{1,\delta_2,\delta_3,\delta_4,\omega,\delta_2\omega,\delta_3\omega,\delta_4\omega\}$ 
we obtain an integral basis of $K$ using the same procedure like in the
proof of Theorem \ref{th3}.

\vspace{0.5cm}

\begin{theorem}
If $K$ admits a power integral basis, then the following divisibility 
conditions must hold:
\label{thx}
\end{theorem}

\noindent
{\scriptsize
\renewcommand{\arraystretch}{2.5}
\begin{tabular}{|c|c|c|c|c|c|c|c|}
\hline n mod 8 & m mod 8 & $1$& $2$ & $3$ & $4$ & $5$ & $6$\\ 
\hline 1,5 & 2,3,6,7 & $n\mid 16m^2\pm1$  & $n\mid 16m\pm1$ & $8m\mid n^4\pm1$    & $8m\mid n^2\pm1$   & $8m\mid 1\pm1$   & $1024m^3\mid m^2(256\pm256)$  \\ 
\hline  1,5 & 1   & $n\mid m^2\pm1$    & $n\mid 4m\pm1$  & $m\mid n^4\pm1$    & $m\mid n^2\pm1$   & $m\mid 1\pm1$   & $16m^3\mid m^2(16\pm16)$  \\ 
\hline 1,5 & 5   & $n\mid m^2\pm1$    & $n\mid 16m\pm1$ & $m\mid n^4\pm1$    & $4m\mid n^2\pm1$   & $m\mid 1\pm1$   & $16m^3\mid m^2(16\pm256)$  \\ 
\hline 2,6 & 1   & $4n\mid m^2\pm1$    & $2n\mid 2m\pm2$  & $m\mid 64n^4\pm1$  & $2m\mid 8n^2\pm2$  & $m\mid 1\pm4$   & $4m^3\mid m^2(16\pm4)$  \\ 
\hline 2,6 & 5   & $4n\mid m^2\pm1$    & $4n\mid 4m\pm4$  & $m\mid 16n^4\pm1$  & $4m\mid 4n^2\pm4$  & $m\mid 1\pm16$  & $16m^3\mid m^2(16\pm16)$  \\
\hline 3,7 & 2,6 & $4n\mid 4m^2\pm4$   & $2n\mid 8m\pm2$  & $16m\mid 4n^4\pm4$   & $8m\mid 2n^2\pm2$  & $16m\mid 4\pm4$   & $256m^3\mid m^2(64\pm64)$  \\
\hline 3,7 & 3,7 & $n\mid m^2\pm16$   & $n\mid 16m\pm1$ & $m\mid n^4\pm16$   & $8m\mid n^2\pm1$   & $m\mid 16\pm1$  & $16m^3\mid m^2(16\pm256)$  \\
\hline 3,7 & 1   & $4n\mid m^2\pm1$    & $n\mid m\pm4$   & $m\mid 16n^4\pm1$  & $4m\mid 4n^2\pm4$  & $m\mid 1\pm16$  & $m^3\mid m^2(16\pm1)$  \\
\hline 3,7 & 5   & $4n\mid m^2\pm1$    & $8n\mid 4m\pm4$  & $m\mid 16n^4\pm1$  & $8m\mid 4n^2\pm4$  & $m\mid 1\pm16$  & $16m^3\mid m^2(16\pm16)$  \\
\hline 2   & 2   & $\displaystyle \frac{n}{2}\mid \frac{m^2}{4}\pm64$ & $\displaystyle \frac{n}{2}\mid 16m\pm1$ & $\displaystyle \frac{m}{2}\mid \frac{n^4}{16}\pm64$& $\displaystyle 8m\mid \frac{n^2}{4}\pm1$ & $\displaystyle \frac{m}{2}\mid 64\pm1$  & $2m^3\mid m^2(4\pm256)$  \\
\hline 2   & 3   & $8n\mid 4m^2\pm4$   & $n\mid 2m\pm8$  & $\displaystyle 8m\mid \frac{n^4}{4}\pm4$  & $\displaystyle 2m\mid \frac{n^2}{2}\pm8$ & $4m\mid 4\pm64$  & $4m^3\mid m^2(64\pm4)$  \\
\hline 2   & 6   & $2n\mid m^2\pm16$   & $n\mid 8m\pm2$  & $\displaystyle 2m\mid \frac{n^4}{4}\pm16$ & $\displaystyle 8m\mid \frac{n^2}{2}\pm2$ & $2m\mid 16\pm4$  & $8m^3\mid m^2(16\pm64)$  \\
\hline 2   & 7   & $8n\mid 4m^2\pm4$   & $\displaystyle \frac{n}{2}\mid m\pm16$  & $\displaystyle m\mid \frac{n^4}{16}\pm4$ & $\displaystyle m\mid \frac{n^2}{4}\pm16$& $4m\mid 4\pm256$ & $m^3\mid m^2(64\pm1)$  \\
\hline 6   & 2   & $2n\mid m^2\pm16$   & $n\mid 8m\pm2$  & $\displaystyle 2m\mid \frac{n^4}{4}\pm16$ & $\displaystyle 8m\mid \frac{n^2}{2}\pm2$ & $2m\mid 16\pm4$  & $8m^3\mid m^2(16\pm64)$  \\
\hline 6   & 3   & $8n\mid 4m^2\pm4$   & $n\mid 2m\pm8$  & $\displaystyle 8m\mid \frac{n^4}{4}\pm4$  & $\displaystyle 2m\mid \frac{n^2}{2}\pm8$ & $4m\mid 4\pm64$  & $4m^3\mid m^2(64\pm4)$  \\
\hline 6   & 6   & $\displaystyle \frac{n}{2}\mid \frac{m^2}{4}\pm64$ & $\displaystyle \frac{n}{2}\mid 16m\pm1$ & $\displaystyle \frac{m}{2}\mid \frac{n^4}{16}\pm64$& $\displaystyle 8m\mid \frac{n^2}{4}\pm1$ & $\displaystyle \frac{m}{2}\mid 64\pm1$  & $2m^3\mid m^2(4\pm256)$  \\
\hline 6   & 7   & $8n\mid 4m^2\pm4$   & $\displaystyle \frac{n}{2}\mid m\pm16$  & $\displaystyle m\mid \frac{n^4}{16}\pm4$ & $\displaystyle m\mid \frac{n^2}{4}\pm16$& $4m\mid 4\pm256$ & $m^3\mid m^2(64\pm1)$ \\
\hline 
\end{tabular} }\\

\vspace{0.5cm}

Note that the formulas like $m|0$ do not yield any 
restriction for $m$, but in all other situations
columns 5 and 6 yield only a few possible values for $m$, and also for $n$ by columns 1 and 2.

In most of the cases column 5 implies that there are only a few
possible values of $m$ and $n$ if $K$ is monogenic. 

\begin{corollary}
If $K$ is monogenic and
$(n\;(\bmod \; 8),m\;(\bmod \; 8))$ is contained in one of the sets
\[
\{1,5\}\times\{5\},\;\;\{2,6\}\times\{1,2,3,5,6,7\},\;\;\{3,7\}\times\{1,3,5,7\},
\]
then
\[
|m|\leq 130 \;\; {\rm and} \;\; |n|\leq 32|m|+32.
\]
\end{corollary}

\vspace{0.5cm}

This yields that the statement is valid in all cases of the above table up to 
lines 1, 2 and 6.

Note that our Corollary above uses only columns 2,5,6 of the above table.

\vspace{0.5cm}

In a recent paper \cite{gr3} the authors considered monogenity of fields
of type $\Q(i,\sqrt[4]{m})$ for squarefree integers
$m\equiv 2,3\;(\bmod \; 4)$. 
Our theorem allows us to extend this 
result to the case $m\equiv 1\;(\bmod \; 4)$, since in these cases
$m|1\pm 16$ must be satisfied:

\begin{corollary}
If $m$ is squarefree integer, 
$|m|\not =1, 3,5,15,17$ then $\Q(i,\sqrt[4]{m})$ not monogenic.
\label{extend}
\end{corollary}

\vspace{0.5cm}

We conjecture that the octic fields $\Q(i,\sqrt[4]{m})$ with 
$m = \pm 3,\pm 5,\pm 15,\pm 17$ are not monogenic, either.
(The fields with $m=\pm 1$ are not octic fields.)

\vspace{0.5cm}

\noindent
{\bf Proof of Theorem \ref{thx}}.\\
Let 
\[
\alpha^{(1)}=\sqrt{n},\;\alpha^{(2)}=-\sqrt{n},\;
\beta^{(1)}=\sqrt[4]{m},\;\beta^{(2)}=i\sqrt[4]{m},\;
\beta^{(3)}=-\sqrt[4]{m},\;\beta^{(4)}=-i\sqrt[4]{m}.
\]
Assume that $\{b_1=1,b_2,\ldots,b_8\}$ is an integral basis of $K$.
The elements of an integral basis are composed of 
$\alpha$ and $\beta$ hence it is unique to denote the conjugate of $b_k$ 
corresponding to $\alpha^{(i)}$ and $\beta^{(j)}$ by $b_k^{(i,j)}$ 
($i=1,2;j=1,2,3$). Set
\[
L^{(i,j)}=L^{(i,j)}(X_1,X_2,\ldots,X_6)=
X_1+b_2^{(i,j)}X_2+\ldots+b_8^{(i,j)}X_8.
\]
Using the $L^{(i,j)}$ we construct the polynomials $F_1,F_2,F_3,F_4,F_5$
in the following way:
\begin{eqnarray*}
F_1&=&\left(L^{(1,1)}-L^{(1,2)}\right)\left(L^{(1,2)}-L^{(1,3)}\right)\left(L^{(1,3)}-L^{(1,4)}\right)\left(L^{(1,4)}-L^{(1,1)}\right)\\
&&\left(L^{(2,1)}-L^{(2,2)}\right)\left(L^{(2,2)}-L^{(2,3)}\right)\left(L^{(2,3)}-L^{(2,4)}\right)\left(L^{(2,4)}-L^{(2,1)}\right),
\end{eqnarray*}
\[
F_2=\left(L^{(1,1)}-L^{(1,3)}\right)\left(L^{(1,2)}-L^{(1,4)}\right)\left(L^{(2,1)}-L^{(2,3)}\right)\left(L^{(2,2)}-L^{(2,4)}\right),
\]
\[
F_3=\left(L^{(1,1)}-L^{(2,1)}\right)\left(L^{(1,2)}-L^{(2,2)}\right)\left(L^{(1,3)}-L^{(2,3)}\right)\left(L^{(1,4)}-L^{(2,4)}\right),
\]
\begin{eqnarray*}
F_4&=&\left(L^{(1,1)}-L^{(2,2)}\right)\left(L^{(1,1)}-L^{(2,4)}\right)\left(L^{(1,2)}-L^{(2,1)}\right)\left(L^{(1,2)}-L^{(2,3)}\right)\\
&&\left(L^{(1,3)}-L^{(2,2)}\right)\left(L^{(1,3)}-L^{(2,4)}\right)\left(L^{(1,4)}-L^{(2,1)}\right)\left(L^{(1,4)}-L^{(2,3)}\right),
\end{eqnarray*}
\[
F_5=\left(L^{(1,1)}-L^{(2,3)}\right)\left(L^{(1,2)}-L^{(2,4)}\right)\left(L^{(1,3)}-L^{(2,1)}\right)\left(L^{(1,4)}-L^{(2,2)}\right).
\]
We find that 
\[
F_i(X_2,\ldots,X_8)=f_i\cdot G_i(X_2,\ldots,X_8) \;\; (i=1,\ldots,5),
\]
where $f_i$ are integers or square roots of integers with $(f_1\ldots f_5)^2=|D_K|$
depending on the parameters only and 
$G_i(X_2,\ldots,X_8) \;\; (i=1,\ldots,5)$ 
are primitive polynomials with integer coefficients.
Then the index form equation corresponding to the basis
$\{1,b_2,\ldots,b_8\}$ as given in Theorem \ref{hhhh} is equivalent to 
\[
G_1(x_2,\ldots,x_8)\ldots G_5(x_2,\ldots,x_8)=\pm 1, \;\; {\rm in}\;\; x_2,\ldots,x_8\in\Z,
\]
that is
\[
G_i(x_2,\ldots,x_8)=\pm 1 \;\; (i=1,\ldots,5),
\]
or equivalently
\[
F_i(x_2,\ldots,x_8)=\pm f_i \;\; (i=1,\ldots,5).
\]
If a power integral basis in $K$ exists then there exist $x_2,\ldots,x_8\in\Z$
satisfying the above equations.
Calculating the explicit form of the $F_i$ in each case of the integral basis
we always have 
\[
n|F_1-F_4, \;\; n|F_2-F_5, \;\;  m|F_3^2-F_4, \;\; m|F_3-F_5, \;\;
m|F_4-F_5^2, \;\; m^3|16F_1-F_2^2,
\]
where the polynomials involved have integer coefficients.
These imply the divisibility conditions of the Theorem.
\mbox{}\\{$\Box$}

\section{Composites of \\quadratic fields and the simplest quartic fields}

Throughout this section we assume that 
\begin{eqnarray}
&&n,m \;\; {\rm are}\;\; {\rm integers}, n\not = 0,1,m\not=0,\pm 3 \;\;{\rm such}\;\;{\rm that} \;\;
n \;\;{\rm is}\;\;{\rm squarefree},\nonumber \\
&&m_0=m^2+16 \;\;{\rm is}\;\;{\rm not}\;\;{\rm divisible}\;\;{\rm by}\;\;{\rm an}\;\;{\rm odd}\;\;{\rm square},\;\;
\gcd(n,m_0)\in\{1,2\}.
\label{f4}
\end{eqnarray}
In this section we consider the composite field $K=\Q(\alpha,\beta)$, where
\begin{center}
\begin{tabular}{rcl}
$\alpha$ &is a root of & $f(x)=x^2-n$,\\
$\beta$  &is a root of & $g(x)=x^4-mx^3-6x^2+mx+1$.\\
\end{tabular}
\end{center}
The fields $M=\Q(\beta)$ are called {\em simplest quartic fields}, 
see \cite{gras}.\\

\begin{theorem}
According to the behaviour of $n$ mod 8, and  $m$ mod 16
an integral basis of $K$ is given by the following table:
\label{th9}
\end{theorem}

\begin{landscape}
{\noindent\scriptsize
\renewcommand{\arraystretch}{2.5}
\begin{tabular}{|c|c|c|}
\hline n mod 8 & m mod 16 & integral basis \\ 
\hline 1,5 & 1,3,5,7,9,11,13,15 & 
$\displaystyle\left\{ 1, \beta, \beta^2, \frac{\beta^3+1}{ 2}, \frac{\alpha+1}{ 2}, \frac{\beta\alpha+\beta}{ 2}, \frac{\beta^2\alpha+\beta^2}{ 2}, \frac{\beta^3\alpha+\beta^3+\alpha+1}{ 4}\right\} $ \\ 
\hline 1,5 & 2,6,10,14 & 
$\displaystyle\left\{1, \beta, \frac{\beta^2+1}{ 2}, \frac{\beta^3+\beta}{ 2}, \frac{\alpha+1}{ 2}, \frac{\beta\alpha+\beta}{ 2}, \frac{\beta^2\alpha+\beta^2+\alpha+1}{ 4}, \frac{\beta^3\alpha+\beta^3+\beta\alpha+\beta}{ 4}\right\} $\\ 
\hline 1,5 & 4,12 & 
$\displaystyle\left\{ 1, \beta, \frac{\beta^2+1}{ 2}, \frac{\beta^3+\beta^2+\beta+1}{ 4}, \frac{\alpha+1}{ 2}, \frac{\beta\alpha+\beta}{ 2}, \frac{\beta^2\alpha+\beta^2+\alpha+1}{ 4}, \frac{\beta^3\alpha+\beta^3+\beta^2\alpha+\beta^2+\beta\alpha+\beta+\alpha+1}{ 8}\right\} $\\ 
\hline 1,5 & 0,8 & 
$\displaystyle\left\{ 1, \beta, \frac{\beta^2+2\beta+3}{ 4}, \frac{\beta^3+3\beta+2}{ 4}, \frac{\alpha+1}{ 2}, \frac{\beta\alpha+\beta}{ 2}, \frac{\beta^2\alpha+\beta^2+2\beta\alpha+2\beta+3\alpha+3}{ 8}, \frac{\beta^3\alpha+\beta^3+3\beta\alpha+3\beta+2\alpha+2}{ 8}
\right\} $\\ 
\hline 2,6 & 1,3,5,7,9,11,13,15 & 
$\displaystyle\left\{1, \beta, \beta^2, \frac{\beta^3+1}{ 2}, \alpha, \alpha\beta, \alpha\beta^2, \frac{\beta^3\alpha+\alpha}{ 2}\right\} $\\ 
\hline 2,6 & 2,6,10,14 & 
$\displaystyle\left\{1, \beta, \frac{\beta^2+1}{ 2}, \frac{\beta^3+\beta}{ 2}, \alpha, \frac{\beta\alpha+\alpha}{ 2}, \frac{\beta^2\alpha+\alpha}{ 2}, \frac{\beta^3\alpha+\beta^2\alpha+\beta\alpha+\alpha}{ 4}\right\} $\\
\hline 2 & 4 & 
$\displaystyle\left\{1, \beta, \frac{\beta^2+1}{ 2}, \frac{\beta^3+\beta^2+\beta+1}{ 4}, \frac{\beta^2+2\beta+2\alpha+3}{ 4}, \frac{\beta^3+\beta^2+2\beta\alpha+\beta+2\alpha+5}{ 8}, \frac{\beta^3+\beta^2\alpha+\beta^2+5\beta+\alpha+5}{ 8}, \frac{\beta^3\alpha+\beta^2\alpha+2\beta^2+\beta\alpha+\alpha+2}{ 8}\right\} $\\
\hline 2 & 12 & 
$\displaystyle\left\{1, \beta, \frac{\beta^2+1}{ 2}, \frac{\beta^3+\beta^2+\beta+1}{ 4}, \frac{\beta^2+2\beta+2\alpha+3}{ 4}, \frac{\beta^3+\beta^2+2\beta\alpha+5\beta+2\alpha+1}{ 8}, \frac{\beta^3+\beta^2\alpha+3\beta^2+5\beta+\alpha+7}{ 8}, \frac{\beta^3\alpha+\beta^2\alpha+2\beta^2+\beta\alpha+\alpha+2}{ 8}\right\} $\\
\hline 2 & 0 & 
$\displaystyle\left\{1, \beta, \frac{\beta^2+2\beta+3}{ 4}, \frac{\beta^3+3\beta+2}{ 4}, \frac{1+\beta+\alpha}{ 2}, \frac{\beta^3+\beta^2+2\beta\alpha+5\beta+2\alpha+5}{ 8}, \frac{\beta^3+\beta^2\alpha+\beta^2+2\beta\alpha+9\beta+3\alpha+5}{ 8}, \frac{\beta^3\alpha+2\beta^3+\beta^2\alpha+9\beta\alpha+14\beta+5\alpha+4}{ 16}\right\} $\\
\hline 2 & 8 & 
$\displaystyle\left\{1, \beta, \frac{\beta^2+2\beta+3}{ 4}, \frac{\beta^3+3\beta+2}{ 4}, \frac{1+\beta+\alpha}{ 2}, \frac{\beta\alpha+\alpha}{ 2}, \frac{\beta^3+\beta^2\alpha+\beta^2+2\beta\alpha+9\beta+3\alpha+5}{ 8}, \frac{\beta^3\alpha+\beta^2\alpha+9\beta\alpha+5\alpha}{ 8}\right\} $\\
\hline 6 & 4 & 
$\displaystyle\left\{1, \beta, \frac{\beta^2+1}{ 2}, \frac{\beta^3+\beta^2+\beta+1}{ 4}, \frac{\beta^2+2\beta+2\alpha+3}{ 4}, \frac{\beta\alpha+\alpha}{ 2}, \frac{\beta^3+\beta^2\alpha+3\beta^2+5\beta+\alpha+7}{ 8}, \frac{\beta^3\alpha+\beta^2\alpha+2\beta^2+\beta\alpha+\alpha+2}{ 8}\right\} $\\
\hline 6 & 12 & 
$\displaystyle\left\{1, \beta, \frac{\beta^2+1}{ 2}, \frac{\beta^3+\beta^2+\beta+1}{ 4}, \frac{\beta^2+2\beta+2\alpha+3}{ 4}, \frac{\beta\alpha+\alpha}{ 2}, \frac{\beta^3+\beta^2\alpha+\beta^2+5\beta+\alpha+5}{ 8}, \frac{\beta^3\alpha+\beta^2\alpha+2\beta^2+\beta\alpha+\alpha+2}{ 8}\right\} $\\
\hline 6 & 0 & 
$\displaystyle\left\{1, \beta, \frac{\beta^2+2\beta+3}{ 4}, \frac{\beta^3+3\beta+2}{ 4}, \frac{1+\beta+\alpha}{ 2}, \frac{\beta\alpha+\alpha}{ 2}, \frac{\beta^3+\beta^2\alpha+\beta^2+2\beta\alpha+9\beta+3\alpha+5}{ 8}, \frac{\beta^3\alpha+\beta^2\alpha+9\beta\alpha+5\alpha}{ 8}\right\} $\\
\hline 6 & 8 & 
$\displaystyle\left\{1, \beta, \frac{\beta^2+2\beta+3}{ 4}, \frac{\beta^3+3\beta+2}{ 4}, \frac{1+\beta+\alpha}{ 2}, \frac{\beta^3+\beta^2+2\beta\alpha+5\beta+2\alpha+5}{ 8}, \frac{\beta^3+\beta^2\alpha+\beta^2+2\beta\alpha+9\beta+3\alpha+5}{ 8}, \frac{\beta^3\alpha+\beta^2\alpha+2\beta^2+13\beta\alpha+4\beta+9\alpha+14}{ 16}\right\} $\\
\hline 3,7 & 1,3,5,7,9,11,13,15 & 
$\displaystyle\left\{1, \beta, \beta^2, \frac{\beta^3+1}{ 2}, \alpha, \alpha\beta, \alpha\beta^2, \frac{\beta^3\alpha+\alpha}{ 2}\right\} $\\
\hline 3,7 & 2,6,10,14 & 
$\displaystyle\left\{1, \beta, \frac{\beta^2+1}{ 2}, \frac{\beta^3+\beta}{ 2}, \frac{\beta^3+\beta^2+3\beta+2\alpha+1}{ 4}, \frac{\beta^3+\beta^2+2\beta\alpha+\beta+3}{ 4}, \frac{\beta^2\alpha+\beta^2+2\beta+\alpha+3}{ 4}, \frac{\beta^3\alpha+\beta^3+\beta\alpha+3\beta+2}{ 4}\right\} $\\
\hline 3,7 & 4,12 & 
$\displaystyle\left\{1, \beta, \frac{\beta^2+1}{ 2}, \frac{\beta^3+\beta^2+\beta+1}{ 4}, \alpha, \frac{\beta\alpha+\beta+\alpha+1}{ 2}, \frac{\beta^2\alpha+\beta^2+\alpha+1}{ 4}, \frac{\beta^3\alpha+\beta^3+\beta^2\alpha+3\beta^2+\beta\alpha+\beta+\alpha+3}{ 8}\right\} $\\
\hline 3,7 & 0,8 & 
$\displaystyle\left\{1, \beta, \frac{\beta^2+2\beta+3}{ 4}, \frac{\beta^3+3\beta+2}{ 4}, \alpha, \frac{\beta\alpha+\beta+\alpha+1}{ 2}, \frac{\beta^2\alpha+2\beta\alpha+3\alpha}{ 4}, \frac{\beta^3\alpha+\beta^3+\beta^2\alpha+\beta^2+5\beta\alpha+5\beta+9\alpha+9}{ 8}\right\} $\\
\hline 
\end{tabular} }\\

\end{landscape}

\vspace{0.5cm}

\noindent
{\bf Proof.}\\
An integral basis $\{1,\omega\}$ and discriminant $D_L$ of $L=\Q(\alpha)$
are given in (\ref{qu}).
According to \cite{kim}, under the conditions of our theorem
an integral basis $\{1,\delta_2,\delta_3,\delta_4\}$
of the simplest quartic fields are given by
\[
\begin{array}{rl}
\left\{1,\beta,\beta^2,\displaystyle\frac{1+\beta^3}{2} \right\}, & \; {\rm if}\;\; \nu_2(m)=0,\\
\left\{1,\beta,\displaystyle\frac{1+\beta^2}{2},\frac{\beta+\beta^3}{2} \right\}, & \; {\rm if}\;\; \nu_2(m)=1,\\
\left\{1,\beta,\displaystyle\frac{1+\beta^2}{2},\frac{1+\beta+\beta^2+\beta^3}{4} \right\}, & \; {\rm if}\;\; \nu_2(m)=2,\\
\left\{1,\beta,\displaystyle\frac{1+2\beta-\beta^2}{4},\frac{1+\beta+\beta^2+\beta^3}{4} \right\}, & \; {\rm if}\;\; \nu_2(m)> 2.\\
\end{array}
\]
The discriminant of $g(x)$ is $4m_0^3$. The conditions on $n$ and $m_0$ 
imply that the discriminant of $K$ is divisible by $n_1^4m_1^2$ where
$m_1=m_0/2^{\nu_2(m_0)}, n_1=n/2^{\nu_2(n)}$ (similarly as in (\ref{reldiscr})).
Starting from the initial basis 
$\{1,\delta_2,\delta_3,\delta_4,\omega,\delta_2\omega,\delta_3\omega,\delta_4\omega\}$ 
we obtain an integral basis of $K$ using the same procedure like in the
proof of Theorem \ref{th3}. Reducing the discriminant of the above initial basis
we only have to deal with 2-factors.
\mbox{}\\{$\Box$}

\vspace{0.5cm}

\begin{theorem}
If $K$ admits a power integral basis, then the following divisibility 
conditions must hold:
\label{thy}
\end{theorem}
\noindent

\hspace*{-1.5cm}
\noindent
{\scriptsize
\renewcommand{\arraystretch}{2.5}
\begin{tabular}{|c|c|c|c|c|c|c|c|}
\hline n mod 8 & m mod 16 & 1& 2 & 3 & 4 & 5 & 6\\ 
\hline 1,5 & 1,3,5,7,11,13,15 & $n\mid m_0^2\pm1$     & $n\mid m_0\pm1$  & $m_0\mid n^4\pm1$   & $m_0\mid n^2\pm1$  & $m_0\mid 1\pm1$    & $m_0^3\mid m_0^2(16\pm1)$ \\ 
\hline 1,5 & 2,6,10,14     & $\displaystyle n\mid \frac{m_0^2}{16}\pm1$  & $n\mid 4m_0\pm1$  & $\displaystyle \frac{m_0}{4}\mid n^4\pm1$   & $m_0\mid n^2\pm1$  & $\displaystyle \frac{m_0}{4}\mid 1\pm1$    & $\displaystyle \frac{m_0^3}{4}\mid m_0^2(1\pm16)$ \\ 
\hline 1,5 & 4,12       & $\displaystyle n\mid \frac{m_0^2}{16}\pm1$  & $n\mid m_0\pm1$   & $\displaystyle \frac{m_0}{2}\mid n^4\pm1$   & $m_0\mid n^2\pm1$  & $\displaystyle \frac{m_0}{2}\mid 1\pm1$    & $\displaystyle \frac{m_0^3}{4}\mid m_0^2(1\pm1)$ \\ 
\hline 1,5 & 0,8       & $\displaystyle n\mid\frac{m_0^2}{256}\pm1$ & $n\mid 4m_0\pm1$  & $\displaystyle \frac{m_0}{16}\mid n^4\pm1$   & $m_0\mid n^2\pm1$  & $\displaystyle \frac{m_0}{16}\mid 1\pm1$    & $\displaystyle \frac{m_0^3}{256}\mid m_0^2\left(\frac{1}{16}\pm16\right)$ \\ 
\hline 2,6,3,7 & 1,3,5,7,9,11,13,15 & $4n\mid m_0^2\pm1$     & $n\mid m_0\pm1$   & $m_0\mid 256n^4\pm1$& $m_0\mid 16n^2\pm1$& $m_0\mid 1\pm1$    & $m_0^3\mid m_0^2(16\pm1)$ \\ 
\hline 2,6 & 2,6,10,14     & $\displaystyle 4n\mid \frac{m_0^2}{16}\pm1$  & $n\mid m_0\pm4$   & $\displaystyle \frac{m_0}{4}\mid 16n^4\pm1$ & $m_0\mid 4n^2\pm4$ & $\displaystyle \frac{m_0}{4}\mid 1\pm16$    & $\displaystyle \frac{m_0^3}{4}\mid m_0^2(1\pm1)$  \\
\hline 2 & 4,12       & $\displaystyle \frac{n}{2}\mid \frac{m_0^2}{1024}\pm64$ & $\displaystyle \frac{n}{2}\mid m_0\pm1$   & $\displaystyle \frac{m_0}{32}\mid \frac{n^4}{16}\pm64$  & $\displaystyle m_0\mid \frac{n^2}{4}\pm1$  & $\displaystyle \frac{m_0}{32}\mid 64\pm1$ & $\displaystyle \frac{m_0^3}{2048}\mid m_0^2\left(\frac{1}{64}\pm1\right)$ \\
\hline 2 & 0       & $\displaystyle 4n\mid \frac{m_0^2}{256}\pm1$ & $\displaystyle \frac{n}{2}\mid \frac{m_0}{16}\pm64$   & $\displaystyle \frac{m_0}{16}\mid \frac{n^4}{16}\pm1$  & $\displaystyle \frac{m_0}{16}\mid \frac{n^2}{4}\pm64$  & $\displaystyle \frac{m_0}{16}\mid 1\pm4096$ & $\displaystyle \frac{m_0^3}{4096}\mid m_0^2\left(\frac{1}{16}\pm\frac{1}{256}\right)$ \\
\hline 2 & 8       & $\displaystyle 4n\mid \frac{m_0^2}{256}\pm1$ & $\displaystyle n\mid \frac{m_0}{4}\pm16$   & $\displaystyle \frac{m_0}{16}\mid n^4\pm1$  & $\displaystyle m_0\mid n^2\pm16$  & $\displaystyle \frac{m_0}{16}\mid 1\pm256$ & $\displaystyle \frac{m_0^3}{4096}\mid m_0^2\left(\frac{1}{16}\pm\frac{1}{16}\right)$ \\
\hline 6 & 4 ,12      & $\displaystyle 2n\mid \frac{m_0^2}{256}\pm16$ & $\displaystyle n\mid \frac{m_0}{2}\pm2$   & $\displaystyle \frac{m_0}{8}\mid \frac{n^4}{4}\pm16$  & $\displaystyle m_0\mid \frac{n^2}{2}\pm2$  & $\displaystyle \frac{m_0}{8}\mid 16\pm4$ & $\displaystyle \frac{m_0^3}{512}\mid m_0^2\left(\frac{1}{16}\pm\frac{1}{4}\right)$ \\
\hline 6 & 0       & $\displaystyle 4n\mid \frac{m_0^2}{256}\pm1$ & $\displaystyle n\mid \frac{m_0}{4}\pm16$   & $\displaystyle \frac{m_0}{16}\mid n^4\pm1$  & $\displaystyle m_0\mid n^2\pm16$  & $\displaystyle \frac{m_0}{16}\mid 1\pm256$ & $\displaystyle \frac{m_0^3}{4096}\mid m_0^2\left(\frac{1}{16}\pm\frac{1}{16}\right)$ \\
\hline 6 & 8       & $\displaystyle 4n\mid \frac{m_0^2}{256}\pm1$ & $\displaystyle \frac{n}{2}\mid \frac{m_0}{16}\pm64$   & $\displaystyle \frac{m_0}{16}\mid \frac{n^4}{16}\pm1$  & $\displaystyle \frac{m_0}{16}\mid \frac{n^2}{4}\pm64$  & $\displaystyle \frac{m_0}{16}\mid 1\pm4096$ & $\displaystyle \frac{m_0^3}{4096}\mid m_0^2\left(\frac{1}{16}\pm\frac{1}{256}\right)$ \\
\hline 3,7 & 2,6,10,14     & $\displaystyle 4n\mid \frac{m_0^2}{16}\pm1$  & $\displaystyle n\mid \frac{m_0}{4}\pm16$& $\displaystyle \frac{m_0}{4}\mid n^4\pm1$   & $\displaystyle \frac{m_0}{4}\mid n^2\pm16$ & $\displaystyle \frac{m_0}{4}\mid 1\pm256$   & $\displaystyle \frac{m_0^3}{64}\mid m_0^2\left(1\pm\frac{1}{16}\right)$ \\
\hline 3,7 & 4,12       & $\displaystyle 4n\mid \frac{m_0^2}{64}\pm4$  & $\displaystyle n\mid \frac{m_0}{2}\pm2$ & $m_0\mid 4n^4\pm4$  & $m_0\mid 2n^2\pm2$ & $m_0\mid 4\pm4$   & $\displaystyle \frac{m_0^3}{8}\mid m_0^2\left(\frac{1}{4}\pm\frac{1}{4}\right)$ \\
\hline 3,7 & 0,8       & $\displaystyle 8n\mid \frac{m_0^2}{256}\pm1$ & $n\mid m_0\pm4$   & $\displaystyle \frac{m_0}{16}\mid 16n^4\pm1$ & $m_0\mid 4n^2\pm4$ & $\displaystyle \frac{m_0}{16}\mid 1\pm16$   & $\displaystyle \frac{m_0^3}{256}\mid m_0^2\left(\frac{1}{16}\pm1\right)$ \\
\hline 
\end{tabular} }

\vspace{0.5cm}

If the right hand term in the divisibility conditions of columns 5 and 6
does not reduce to 0, there remain only a few possible values for $m_0$. 

\begin{corollary}
If $(n\mod 8,\; m\mod 16)\not \in \{1,3,5,7\}\times\{4,12\}$ and $K$ is monogenic, then
\[
|m|\leq 64 \;\;{\rm and}\;\; |n|\leq 4m^2+192.
\]
\end{corollary}

\vspace{0.5cm}

Note that our Corollary uses only columns 2,5,6 of the above table.

\vspace{0.5cm}

\noindent
{\bf Proof of Theorem \ref{thy}.}\\
Let 
\[
\alpha^{(1)}=\sqrt{n},\;\alpha^{(2)}=-\sqrt{n},\;
\beta^{(1)}=\beta,\;\beta^{(2)}=\frac{\beta-1}{\beta+1},\;
\beta^{(3)}=\frac{-1}{\beta},\;\beta^{(4)}=\frac{-\beta-1}{\beta-1}.
\]
Assume that $\{b_1=1,b_2,\ldots,b_8\}$ is the integral basis of $K$
as given in Theorem \ref{th9}.
The elements of the integral basis are composed of 
$\alpha$ and $\beta$ hence it is unique to denote the conjugate of $b_k$ 
corresponding to $\alpha^{(i)}$ and $\beta^{(j)}$ by $b_k^{(i,j)}$ 
($i=1,2;j=1,2,3,4$). We construct the same $L^{(i,j)}$ and use
the analogous factorization of the index form into
polynomials $F_1,F_2,F_3,F_4,F_5$ like in the proof of Theorem \ref{thx}.
Similarly, the properties
\[
n|F_1-F_4, \;\; n|F_2-F_5, \;\;  m_0|F_3^2-F_4, \;\; m_0|F_3-F_5, \;\;
m_0|F_4-F_5^2, \;\; m_0^3|16F_1-F_2^2
\]
(where the polynomials involved have integer coefficients)
imply the assertions of our Theorem.
\mbox{}\\{$\Box$}

\vspace{0.5cm}

\section{The composite fields $\Q(i\sqrt{3},\sqrt[6]{m})$}

M.-L. Chang \cite{chang} considered the normal closure of pure cubic fields,
the number fields $\Q(\omega,\sqrt[3]{m})$, where $\omega=e^{2\pi i/3}$
and $m$ is not a complete cube.
He showed that these fields are not monogenic, except for $m=2$.

In a recent paper \cite{gr3} the authors studied the analogous problem 
in the fields $\Q(i,\sqrt[4]{m})$ and proved that if 
$m$ is a square-free integer, $m\equiv 2,3 \; (\bmod \; {4})$, then
the field $\Q(i,\sqrt[4]{m})$ is not monogenic. This assertion was extended
in Corollary \ref{extend} of the present paper.\\

Using the technics developed in our paper we shall now continue the series 
of the above
mentioned results and consider the fields $K=\Q(\omega,\sqrt[6]{m})$
where $\omega=e^{2\pi i/3}$. 

\noindent
Throughout this section we assume that
\begin{equation}
m \;\;{\rm is}\;\; {\rm a}\;\;{\rm square}\;\;{\rm free}\;\;{\rm integer}, m\not =0,\pm 1, -3.
\label{f5}
\end{equation}

\begin{theorem}
According to the behaviour of $m$ mod 36
an integral basis of $K$ is given by the following table
\label{th12}
\end{theorem}

\begin{landscape}
\vspace*{-1cm}
{\tiny
\renewcommand{\arraystretch}{2.5}
\begin{tabular}{|c|c|}
\hline m mod 36 & integral basis\\ 
\hline 1 &  
$\displaystyle\left\{ 1, \beta, \beta^2, \frac{\beta^3+1}{ 2}, \frac{\beta^4+4\beta^2+3\beta+4}{ 6}, \frac{\beta^5+\beta^3+3\beta^2+4\beta+3}{ 6}, \omega, \omega\beta, \frac{\beta^2\omega+\beta^2+2\omega+2}{ 3}, \frac{\beta^3\omega+\beta^3+2\beta\omega+2\beta+3\omega+3}{ 6}, \frac{\beta^4\omega+4\beta^2\omega+3\beta\omega+4\omega}{ 6}, \frac{\beta^5\omega+\beta^3\omega+3\beta^2\omega+4\beta\omega+3\omega}{ 6}\right\} $ \\ 
\hline 2 &  
$\displaystyle\left\{1, \beta, \beta^2, \beta^3, \beta^4, \beta^5, \omega, \omega\beta, \omega\beta^2, \omega\beta^3, \frac{\beta^4\omega+\beta^4+2\beta^2\omega+2\beta^2+\omega+1}{ 3}, \frac{\beta^5\omega+\beta^5+2\beta^3\omega+2\beta^3+\beta\omega+\beta}{ 3}\right\} $\\ 
\hline 3 &  
$\displaystyle\left\{ 1, \beta, \beta^2, \beta^3, \beta^4, \beta^5, \omega, \omega\beta, \omega\beta^2, \frac{\beta^3\omega+\beta^3}{ 3}, \frac{\beta^4\omega+\beta^4}{ 3}, \frac{\beta^5\omega+\beta^5}{ 3}\right\} $\\ 
\hline 5 &
$\displaystyle\left\{ 1, \beta, \beta^2, \frac{\beta^3+1}{ 2}, \frac{\beta^4+\beta}{ 2}, \frac{\beta^5+\beta^2}{ 2}, \omega, \omega\beta, \omega\beta^2, \frac{\beta^3\omega+\omega}{ 2}, \frac{\beta^4\omega+\beta^4+2\beta^2\omega+2\beta^2+3\beta\omega+3\beta+4\omega+4}{ 6}, \frac{\beta^5\omega+\beta^5+2\beta^4\omega+2\beta^4+2\beta^3\omega+2\beta^3+\beta^2\omega+\beta^2+4\beta\omega+4\beta+2\omega+2}{ 6}\right\} $\\ 
\hline 6 & 
$\displaystyle\left\{1, \beta, \beta^2, \beta^3, \beta^4, \beta^5, \omega, \omega\beta, \omega\beta^2, \frac{\beta^3\omega+\beta^3}{ 3}, \frac{\beta^4\omega+\beta^4}{ 3}, \frac{\beta^5\omega+\beta^5}{ 3}\right\} $\\ 
\hline 7 &
$\displaystyle\left\{
1, \beta, \beta^2, \beta^3, \beta^4, \beta^5, \omega, \omega\beta, \omega\beta^2, \omega\beta^3, \frac{\beta^4\omega+\beta^4+\beta^2\omega+\beta^2+\omega+1}{ 3}, \frac{\beta^5\omega+\beta^5+\beta^3\omega+\beta^3+\beta\omega+\beta}{ 3} \right\} $\\
\hline 10 & 
$\displaystyle\left\{
1, \beta, \beta^2, \beta^3, \frac{\beta^4+\beta^2+1}{ 3}, \frac{\beta^5+\beta^3+\beta}{ 3}, \omega, \omega\beta, \frac{\beta^2\omega+\beta^2+2\omega+2}{ 3}, \frac{\beta^3\omega+\beta^3+2\beta\omega+2\beta}{ 3}, \frac{\beta^4\omega+\beta^2\omega+\omega}{ 3}, \frac{\beta^5\omega+\beta^3\omega+\beta\omega}{ 3}
\right\} $\\
\hline 11 &  
$\displaystyle\left\{1, \beta, \beta^2, \beta^3, \beta^4, \beta^5, \omega, \omega\beta, \omega\beta^2, \omega\beta^3, \frac{\beta^4\omega+\beta^4+2\beta^2\omega+2\beta^2+\omega+1}{ 3}, \frac{\beta^5\omega+\beta^5+2\beta^3\omega+2\beta^3+\beta\omega+\beta}{ 3}
\right\} $\\
\hline 13 & 
$\displaystyle\left\{1, \beta, \beta^2, \frac{\beta^3+1}{ 2}, \frac{\beta^4+\beta}{ 2}, \frac{\beta^5+\beta^2}{ 2}, \omega, \omega\beta, \omega\beta^2, \frac{\beta^3\omega+\omega}{ 2}, \frac{\beta^4\omega+\beta^4+4\beta^2\omega+4\beta^2+3\beta\omega+3\beta+4\omega+4}{ 6}, \frac{\beta^5\omega+\beta^5+\beta^4\omega+\beta^4+\beta^3\omega+\beta^3+\beta^2\omega+\beta^2+\beta\omega+\beta+\omega+1}{ 6}\right\} $\\
\hline 14 & 
$\displaystyle\left\{
1, \beta, \beta^2, \beta^3, \beta^4, \beta^5, \omega, \omega\beta, \omega\beta^2, \omega\beta^3, \frac{\beta^4\omega+\beta^4+2\beta^2\omega+2\beta^2+\omega+1}{ 3}, \frac{\beta^5\omega+\beta^5+2\beta^3\omega+2\beta^3+\beta\omega+\beta}{ 3}\right\} $\\
\hline 15 & 
$\displaystyle\left\{
1, \beta, \beta^2, \beta^3, \beta^4, \beta^5, \omega, \omega\beta, \omega\beta^2, \frac{\beta^3\omega+\beta^3}{ 3}, \frac{\beta^4\omega+\beta^4}{ 3}, \frac{\beta^5\omega+\beta^5}{ 3}
\right\} $\\
\hline 17 & 
$\displaystyle\left\{1, \beta, \beta^2, \frac{\beta^3+1}{ 2}, \frac{\beta^4+2\beta^2+3\beta+4}{ 6}, \frac{\beta^5+2\beta^3+3\beta^2+4\beta}{ 6}, \omega, \omega\beta, \frac{\beta^2\omega+\beta^2+\omega+1}{ 3}, \frac{\beta^3\omega+\beta^3+4\beta\omega+4\beta+3\omega+3}{ 6}, \frac{\beta^4\omega+2\beta^2\omega+3\beta\omega+4\omega}{ 6}, \frac{\beta^5\omega+2\beta^3\omega+3\beta^2\omega+4\beta\omega}{ 6}
\right\} $\\
\hline 19 &
$\displaystyle\left\{1, \beta, \beta^2, \beta^3, \frac{\beta^4+\beta^2+1}{ 3}, \frac{\beta^5+\beta^3+\beta}{ 3}, \omega, \omega\beta, \frac{\beta^2\omega+\beta^2+2\omega+2}{ 3}, \frac{\beta^3\omega+\beta^3+2\beta\omega+2\beta}{ 3}, \frac{\beta^4\omega+\beta^2\omega+\omega}{ 3}, \frac{\beta^5\omega+\beta^3\omega+\beta\omega}{ 3}\right\} $\\
\hline 21 &  
$\displaystyle\left\{1, \beta, \beta^2, \frac{\beta^3+1}{ 2}, \frac{\beta^4+\beta}{ 2}, \frac{\beta^5+\beta^2}{ 2}, \omega, \omega\beta, \omega\beta^2, \frac{\beta^3\omega+\beta^3+3\omega+3}{ 6}, \frac{\beta^4\omega+\beta^4+3\beta\omega+3\beta}{ 6}, \frac{\beta^5\omega+\beta^5+3\beta^2\omega+3\beta^2}{ 6}
\right\} $\\
\hline 22 &  
$\displaystyle\left\{1, \beta, \beta^2, \beta^3, \beta^4, \beta^5, \omega, \omega\beta, \omega\beta^2, \omega\beta^3, \frac{\beta^4\omega+\beta^4+\beta^2\omega+\beta^2+\omega+1}{ 3}, \frac{\beta^5\omega+\beta^5+\beta^3\omega+\beta^3+\beta\omega+\beta}{ 3}\right\} $\\
\hline 23 &  
$\displaystyle\left\{
1, \beta, \beta^2, \beta^3, \beta^4, \beta^5, \omega, \omega\beta, \omega\beta^2, \omega\beta^3, \frac{\beta^4\omega+\beta^4+2\beta^2\omega+2\beta^2+\omega+1}{ 3}, \frac{\beta^5\omega+\beta^5+2\beta^3\omega+2\beta^3+\beta\omega+\beta}{ 3}
\right\} $\\
\hline 25 &  
$\displaystyle\left\{1, \beta, \beta^2, \frac{\beta^3+1}{ 2}, \frac{\beta^4+\beta}{ 2}, \frac{\beta^5+\beta^2}{ 2}, \omega, \omega\beta, \omega\beta^2, \frac{\beta^3\omega+\omega}{ 2}, \frac{\beta^4\omega+\beta^4+4\beta^2\omega+4\beta^2+3\beta\omega+3\beta+4\omega+4}{ 6}, \frac{\beta^5\omega+\beta^5+\beta^4\omega+\beta^4+\beta^3\omega+\beta^3+\beta^2\omega+\beta^2+\beta\omega+\beta+\omega+1}{ 6}
\right\} $\\
\hline 26 &  
$\displaystyle\left\{1, \beta, \beta^2, \beta^3, \frac{\beta^4+2\beta^2+1}{ 3}, \frac{\beta^5+2\beta^3+\beta}{ 3}, \omega, \omega\beta, \frac{\beta^2\omega+\beta^2+\omega+1}{ 3}, \frac{\beta^3\omega+\beta^3+\beta\omega+\beta}{ 3}, \frac{\beta^4\omega+2\beta^2\omega+\omega}{ 3}, \frac{\beta^5\omega+2\beta^3\omega+\beta\omega}{ 3}\right\} $\\
\hline 29 &  
$\displaystyle\left\{1, \beta, \beta^2, \frac{\beta^3+1}{ 2}, \frac{\beta^4+\beta}{ 2}, \frac{\beta^5+\beta^2}{ 2}, \omega, \omega\beta, \frac{\beta^4\omega+\beta^4+2\beta^2\omega+2\beta^2+3\beta\omega+3\beta+4\omega+4}{ 6}, \frac{\beta^5\omega+\beta^5+2\beta^4\omega+2\beta^4+2\beta^3\omega+2\beta^3+\beta^2\omega+\beta^2+4\beta\omega+4\beta+2\omega+2}{ 6}, \frac{\beta^4\omega+\beta\omega}{ 2}, \frac{\beta^5\omega+\beta^2\omega}{ 2}
\right\} $\\
\hline 30 &  
$\displaystyle\left\{
1, \beta, \beta^2, \beta^3, \beta^4, \beta^5, \omega, \omega\beta, \omega\beta^2, \frac{\beta^3\omega+\beta^3}{ 3}, \frac{\beta^4\omega+\beta^4}{ 3}, \frac{\beta^5\omega+\beta^5}{ 3}\right\} $\\
\hline 31 &  
$\displaystyle\left\{1, \beta, \beta^2, \beta^3, \beta^4, \beta^5, \omega, \omega\beta, \omega\beta^2, \omega\beta^3, \frac{\beta^4\omega+\beta^4+\beta^2\omega+\beta^2+\omega+1}{ 3}, \frac{\beta^5\omega+\beta^5+\beta^3\omega+\beta^3+\beta\omega+\beta}{ 3}\right\} $\\
\hline 33 &  
$\displaystyle\left\{1, \beta, \beta^2, \frac{\beta^3+1}{ 2}, \frac{\beta^4+\beta}{ 2}, \frac{\beta^5+\beta^2}{ 2}, \omega, \omega\beta, \omega\beta^2, \frac{\beta^3\omega+\beta^3+3\omega+3}{ 6}, \frac{\beta^4\omega+\beta^4+3\beta\omega+3\beta}{ 6}, \frac{\beta^5\omega+\beta^5+3\beta^2\omega+3\beta^2}{ 6}
\right\} $\\
\hline 34 &  
$\displaystyle\left\{1, \beta, \beta^2, \beta^3, \beta^4, \beta^5, \omega, \omega\beta, \omega\beta^2, \omega\beta^3, \frac{\beta^4\omega+\beta^4+\beta^2\omega+\beta^2+\omega+1}{ 3}, \frac{\beta^5\omega+\beta^5+\beta^3\omega+\beta^3+\beta\omega+\beta}{ 3}\right\} $\\
\hline 35 &  
$\displaystyle\left\{1, \beta, \beta^2, \beta^3, \frac{\beta^4+2\beta^2+1}{ 3}, \frac{\beta^5+2\beta^3+\beta}{ 3}, \omega, \omega\beta, \frac{\beta^2\omega+\beta^2+\omega+1}{ 3}, \frac{\beta^3\omega+\beta^3+\beta\omega+\beta}{ 3}, \frac{\beta^4\omega+2\beta^2\omega+\omega}{ 3}, \frac{\beta^5\omega+2\beta^3\omega+\beta\omega}{ 3}\right\} $\\
\hline 
\end{tabular} }
\end{landscape}

\vspace{0.5cm}

\noindent
{\bf Proof.}\\
An integral basis of $L=\Q(\omega)$ is $\{1,\omega\}$ and $D_L=-3$.
Set $\beta=\sqrt[6]{m}$.
An integral basis of $M=\Q(\sqrt[6]{m})$ is given in \cite{gr2}:
\[
\begin{array}{rl}
\{1,\beta,\beta^2,\beta^3,\beta^4,\beta^5\}, 
& \;\;{\rm if}\;\; m\equiv 2,3,6,7,11,14,15,22,23,30,31,34 \;(\bmod{\; 36}),\\
\left\{1,\beta,\beta^2,\frac{1+\beta^3}{2},\frac{\beta+\beta^4}{2},
\frac{\beta^2+\beta^5}{2}\right\},
& \;\;{\rm if}\;\; m\equiv 5,13,21,25,29,33 \;(\bmod{\; 36}),\\
\left\{1,\beta,\beta^2, \beta^3, \frac{1+\beta^2+\beta^4}{3},
\frac{\beta+\beta^3+\beta^5}{3}\right\},
& \;\;{\rm if}\;\; m\equiv 10,19 \;(\bmod{\; 36}),\\
\left\{1,\beta,\beta^2, \beta^3,\frac{1+2\beta^2+\beta^4}{3},\frac{\beta+2\beta^3+\beta^5}{3}\right\},
& \;\;{\rm if}\;\; m\equiv 26,35 \;(\bmod{\; 36}),\\
\left\{1,\beta,\beta^2,\frac{1+\beta^3}{2}, \frac{4+3\beta+2\beta^2+\beta^4}{6},
\frac{4\beta+3\beta^2+2\beta^3+\beta^5}{6}\right\},
& \;\;{\rm if}\;\; m\equiv 17 \;(\bmod{\; 36}),\\
\left\{1,\beta,\beta^2,\frac{1+\beta^3}{2},\frac{4+3\beta+4\beta^2+\beta^4}{6},
\frac{3+4\beta+3\beta^2+\beta^3+\beta^5}{6}\right\},
& \;\;{\rm if}\;\; m\equiv 1 \;(\bmod{\; 36}).\\
\end{array}
\]
The discriminant of $g(x)=x^6-m$ is $2^63^6m^5$. Let 
$m_1=m/(2^{\nu_2(m)}3^{\nu_3(m)})$. Similarly as in the previous proofs
$m_1^{10}$ divides $D_K$. 
If $\{1,\delta_2,\ldots,\delta_6\}$ is an integral basis of $M$, then we
start from the initial basis
$\{1,\delta_2,\ldots,\delta_6,\omega,\delta_2\omega,\ldots,\delta_6\omega\}$
of $K$ and consider possible reductions of the discriminant of this basis
by 2-factors and 3-factors.
\mbox{}\\{$\Box$}

\vspace{0.5cm}

\begin{theorem}
If $K$ admits a power integral basis, then the following divisibility 
conditions must hold:
\label{thz}
\end{theorem}

{\scriptsize
\renewcommand{\arraystretch}{2.5}
\begin{tabular}{|c|c|c|c|c|c|}
\hline m mod 36 & 1 & 2 & 3 & 4 & 5 \\ \hline&&&& 6 & 7  \\ \hline
\hline 1,17 &  $m\mid 6561\pm1 $ & $m\mid 6561\pm9 $ & $m\mid 81\pm1 $  &  $m\mid 9\pm 1$ & $m^3\mid m^2(729 \pm 1) $ \\ \hline&&&& $m^3\mid m^2(4096 \pm 1) $ & $m^3\mid m^2(4096 \pm 729) $ \\ \hline
\hline 2,7,11,14,22,23,31,34 &  $m\mid 6561\pm1 $ & $m\mid 6561\pm9 $ & $m\mid 81\pm1 $  &  $m\mid 9\pm 1$ & $m^3\mid m^2(729 \pm 81) $ \\ \hline&&&& $m^3\mid m^2(4096 \pm 4096) $ & $m^3\mid m^2(331776 \pm 2985984) $  \\ \hline
\hline 3,6,15,30 &  $\displaystyle \frac{m}{3}\mid 729\pm9 $ & $\displaystyle \frac{m}{3}\mid 729\pm1 $ & $\displaystyle \frac{m}{3}\mid 27\pm3 $ & $\displaystyle \frac{m}{3}\mid 9\pm 1$  & $\displaystyle \frac{m^3}{27}\mid m^2(81 \pm 729) $ \\ \hline&&&& $ \displaystyle \frac{m^3}{27}\mid \left(\frac{4096m^2}{9} \pm \frac{4096m^2}{9}\right) $ & $\displaystyle \frac{m^3}{27}\mid m^2(2985984 \pm 331776) $    \\ \hline
\hline 5,13,25,29 &  $m\mid 6561\pm1 $ & $m\mid 6561\pm9 $ & $m\mid 81\pm1 $  &  $m\mid 9\pm 1$ & $m^3\mid m^2(729 \pm 81) $ \\ \hline&&&& $m^3\mid m^2(4096 \pm 1) $ & $m^3\mid m^2(331776 \pm 729) $   \\ \hline
\hline 10,19,26,35 &  $m\mid 6561\pm1 $ & $m\mid 6561\pm9 $ & $m\mid 81\pm1 $  &  $m\mid 9\pm 1$ & $m^3\mid m^2(729 \pm 1) $ \\ \hline &&&& $m^3\mid m^2(4096 \pm 4096) $ & $m^3\mid m^2(4096 \pm 2985984) $   \\ \hline
\hline 21,33 &  $\displaystyle \frac{m}{3}\mid 729\pm9 $ & $\displaystyle \frac{m}{3}\mid 729\pm1 $ & $\displaystyle \frac{m}{3}\mid 27\pm3 $ & $\displaystyle \frac{m}{3}\mid 9\pm 1$   & $\displaystyle \frac{m^3}{27}\mid m^2(81 \pm 729) $ \\ \hline&&&& $\displaystyle \frac{m^3}{27}\mid \left(\frac{4096m^2}{9} \pm \frac{4096m^2}{9}\right) $ & $\displaystyle \frac{m^3}{27}\mid m^2(2985984 \pm 81) $    \\ 
\hline 
\end{tabular} 
}

\vspace{0.5cm}

Obviously the conditions can only be satisfied by a few values of $m$,
cf. eg. column 4 of the table. 
Testing the possible values we find:

\begin{corollary}
If $|m|\not = 2,3,5,6,10,15,30$ then $K$ is not monogenic.
\end{corollary}

We conjecture that the field $K$ is
not monogenic for the above values of $m$, either.
Note that our Corollary uses only column 4 of the above table.

\vspace{0.5cm}

\noindent
{\bf Proof of Theorem \ref{thz}}.\\
Let 
\[
\omega^{(1)}=\frac{1+i\sqrt{3}}{2}=\omega,\;\; \omega^{(2)}=\frac{1-i\sqrt{3}}{2},
\]
and
\[
\beta^{(j)}=\sqrt[6]{m}\cdot \omega^{j-1},\;\; j=1,\ldots,6.
\]
Assume that $\{b_1=1,b_2,\ldots,b_{12}\}$ is the integral basis of $K$
as given in Theorem \ref{th12}.
The elements of the integral basis are composed of 
$\omega$ and $\beta$ hence it is unique to denote the conjugate of $b_k$ 
corresponding to $\omega^{(i)}$ and $\beta^{(j)}$ by $b_k^{(i,j)}$ 
($i=1,2;j=1,\ldots,6$). Set
\[
L^{(i,j)}=L^{(i,j)}(X_1,\ldots,X_{12})=
X_1+b_2^{(i,j)}X_2+\ldots+b_{12}^{(i,j)}X_{12}.
\]
Using the $L^{(i,j)}$ we construct the polynomials $F_1,\ldots,F_7$
in the following way:
\begin{eqnarray*}
F_1&=&\left(L^{(1,1)}-L^{(1,2)}\right)\left(L^{(1,2)}-L^{(1,3)}\right)\left(L^{(1,3)}-L^{(1,4)}\right)
\left(L^{(1,4)}-L^{(1,5)}\right)\left(L^{(1,5)}-L^{(1,6)}\right)\left(L^{(1,6)}-L^{(1,1)}\right) \\
&\times&\left(L^{(2,1)}-L^{(2,2)}\right)\left(L^{(2,2)}-L^{(2,3)}\right)\left(L^{(2,3)}-L^{(2,4)}\right)
\left(L^{(2,4)}-L^{(2,5)}\right)\left(L^{(2,5)}-L^{(2,6)}\right)\left(L^{(2,6)}-L^{(2,1)}\right),\\
F_2&=&\left(L^{(1,1)}-L^{(1,3)}\right)\left(L^{(1,2)}-L^{(1,4)}\right)\left(L^{(1,3)}-L^{(1,5)}\right)
\left(L^{(1,4)}-L^{(1,6)}\right)\left(L^{(1,5)}-L^{(1,1)}\right)\left(L^{(1,6)}-L^{(1,2)}\right)\\
&\times&\left(L^{(2,1)}-L^{(2,3)}\right)\left(L^{(2,2)}-L^{(2,4)}\right)\left(L^{(2,3)}-L^{(2,5)}\right)
\left(L^{(2,4)}-L^{(2,6)}\right)\left(L^{(2,5)}-L^{(2,1)}\right)\left(L^{(2,6)}-L^{(2,2)}\right),\\
F_3&=&\left(L^{(1,1)}-L^{(1,4)}\right)\left(L^{(1,2)}-L^{(1,5)}\right)\left(L^{(1,3)}-L^{(1,6)}\right)
\left(L^{(2,1)}-L^{(2,4)}\right)\left(L^{(2,2)}-L^{(2,5)}\right)\left(L^{(2,3)}-L^{(2,6)}\right),\\
F_4&=&\left(L^{(1,1)}-L^{(2,1)}\right)\left(L^{(1,2)}-L^{(2,2)}\right)\left(L^{(1,3)}-L^{(2,3)}\right)
\left(L^{(1,4)}-L^{(2,4)}\right)\left(L^{(1,5)}-L^{(2,5)}\right)\left(L^{(1,6)}-L^{(2,6)}\right),\\
F_5&=&\left(L^{(1,1)}-L^{(2,2)}\right)\left(L^{(1,2)}-L^{(2,3)}\right)\left(L^{(1,3)}-L^{(2,4)}\right)
\left(L^{(1,4)}-L^{(2,5)}\right)\left(L^{(1,5)}-L^{(2,6)}\right)\left(L^{(1,6)}-L^{(2,1)}\right)\\
&\times&\left(L^{(1,1)}-L^{(2,6)}\right)\left(L^{(1,2)}-L^{(2,1)}\right)\left(L^{(1,3)}-L^{(2,2)}\right)
\left(L^{(1,4)}-L^{(2,3)}\right)\left(L^{(1,5)}-L^{(2,4)}\right)\left(L^{(1,6)}-L^{(2,5)}\right),\\
F_6&=&\left(L^{(1,1)}-L^{(2,3)}\right)\left(L^{(1,2)}-L^{(2,4)}\right)\left(L^{(1,3)}-L^{(2,5)}\right)
\left(L^{(1,4)}-L^{(2,6)}\right)\left(L^{(1,5)}-L^{(2,1)}\right)\left(L^{(1,6)}-L^{(2,2)}\right)\\
&\times&\left(L^{(1,1)}-L^{(2,5)}\right)\left(L^{(1,2)}-L^{(2,6)}\right)\left(L^{(1,3)}-L^{(2,1)}\right)
\left(L^{(1,4)}-L^{(2,2)}\right)\left(L^{(1,5)}-L^{(2,3)}\right)\left(L^{(1,6)}-L^{(2,4)}\right),\\
F_7&=&\left(L^{(1,1)}-L^{(2,4)}\right)\left(L^{(1,2)}-L^{(2,5)}\right)\left(L^{(1,3)}-L^{(2,6)}\right)
\left(L^{(1,4)}-L^{(2,1)}\right)\left(L^{(1,5)}-L^{(2,2)}\right)\left(L^{(1,6)}-L^{(2,3)}\right).
\end{eqnarray*}
We find that
\[
F_i(X_2,\ldots,X_{12})=f_i\cdot G_i(X_2,\ldots,X_{12}), \;\; (i=1,\ldots,7),
\]
where $f_i$ are integers or square roots of integers 
depending on $m$ with $(f_1\ldots f_7)^2=|D_K|$
and $G_i$ are primitive polynomials with integer coefficients.
Then the index form equation corresponding to the basis
$\{1,b_2,\ldots,b_{12}\}$ is equivalent to 
\[
G_1(x_2,\ldots,x_{12})\ldots G_7(x_2,\ldots,x_{12})=\pm 1 \;\; {\rm in} \;\; x_2,\ldots,x_{12}\in\Z,
\]
that is
\[
G_i(x_2,\ldots,x_{12})=\pm 1 \;\; (i=1,\ldots,7),
\]
or equivalently
\[
F_i(x_2,\ldots,x_{12})=\pm f_i \;\; (i=1,\ldots,7).
\]
If a power integral basis in $K$ exists then there exist $x_2,\ldots,x_{12}\in\Z$
satisfying the above equations.
Direct calculation shows that
\[
m\mid F_4^2-F_5,\;\; m\mid F_4^2-F_6,\;\;m\mid F_4-F_7,\;\;m\mid F_5-F_6,\;\;
\]
\[
m^3\mid 729F_1-F_2,\;m^3\mid 4096F_1-F_3^2,\;\;m^3\mid 4096F_2-729F_3^2,
\]
(where the polynomials involved have integer coefficients)
which imply the assertion of our Theorem \ref{thz}.
\mbox{}\\{$\Box$}

\section{Computational remarks}

Calculating an integral basis by performing the reduction
by 2-factors and 3-factors gave a huge number of possible
cases and types of the integral basis which are shown in our tables.

We had to calculate the factors of the index form
in all cases of the integral basis explicitly.  
In a number field of degree $n$ the index form has $n-1$
variables and degree $n(n-1)/2$. For $n=6,8,12$ this degree 
is $15,28,66$, respectively. Therefore it was impossible to
calculate the complete index form and then to factorize it.
We had to combine suitable linear factors such that their product is
invariant under some subgroup of the Galois group of $K/\Q$.
Then we explicitely calculated the products of these linear factors.
These products having integer coefficients are the factors of the index form.
Finally we had to find explicitely the divisors of these factors.

This procedure was supported by our experience but it was still a 
huge computation, in a large number of cases. 
The calculation in one case of the integral basis took only some seconds
except in the last section with fields of degree 12 when
it took about 60 minutes per case. 

All calculations were performed in Maple \cite{maple}.
A very careful and efficient handling of the formulas
was absolutely necessary to manage the calculations.


\begin{thebibliography}{10}

\normalsize
\baselineskip=17pt

\bibitem{chang}
Mu-Ling Chang, {\em Non-monogenity in a family of sextic fields},
J. Number Theory, {\bf 97}(2002), 252--268.


\bibitem{maple}
B.W.Char, K.O.Geddes, G.H.Gonnet, M.B.Monagan, S.M.Watt (eds.)
MAPLE, Reference Manual, Watcom Publications, Waterloo, Canada, 1988.

\bibitem{cohen}
H,Cohen,
A course in computational algebraic number theory,
Graduate Texts in Mathematics. 138, Springer, Berlin, 1993.

\bibitem{cook}
J.P.Cook, {\em Computing Integral Bases},
http://math.ou.edu/~jcook/LaTeX/integralbases.pdf .


\bibitem{comp1}
I.Ga\'al, {\em Power integral bases in composites of number fields}, 
Canad. Math. Bulletin {\bf 41}(1998), 158--165.

\bibitem{gaal}
I.Ga\'al, 
Diophantine equations and power integral bases,
Boston, Birkh\"auser, 2002.

\bibitem{comp3}
I.Ga\'al and P.Olajos
{\em Recent results on power integral bases of composite fields},
Acta Acad. Paedagog. Agriensis, Sect. Mat. (N.S.) {\bf 30}(2003), 45-54.

\bibitem{comp2}
I.Ga\'al, P.Olajos and M.Pohst,
{\em Power integer bases in orders of composite fields}, 
Experimental Math. {\bf 11}(2002), No.1., 87--90.

\bibitem{gr2}
I.Ga\'al and L.Remete,
{\em Integral bases and monogenity of pure fields},
J. Number Theory {\bf 173}(2017), 129-146.

\bibitem{gr3}
I.Ga\'al and L.Remete,
{\em Non-monogenity in a family of octic fields},
Rocky Mountain J. Math., to appear.

\bibitem{gr4}
I.Ga\'al and L.Remete,
{\em Integral bases and monogenity of the simplest sextic fields},
submitted.

\bibitem{gras}
M. N. Gras, {\em Table numerique du nombre de
classes et des unites des extensions cycliques reelles de
degr\'e 4 in Q}, Publ. Math. Fac. Sci. Besan¸con (1977-
1978) fasc. 2

\bibitem{kim}
H. K. Kim and J. S. Kim,
{\em Computation of the Different of the Simplest Quartic Fields},
Manuscript, 2003.

\bibitem{nark}
W.Narkiewicz, {\em Elementary and Analytic Theory of Algebraic
Numbers}, Third Edition, Springer, 1990.


\bibitem{comp4}
P.Olajos,
{\em Power integral bases in orders of composite fields. II},
Ann. Univ. Sci. Budap. Rolando Eötvös, Sect. Math. {\bf 46}(2003), 35-41.

\bibitem{shanks}
D.Shanks, {\em The simplest cubic fields}, Math. Comput., {\bf 28}(1974), 
1137--1152.

\bibitem{stewart}
I.Stewart and D.Tall,
Algebraic number theory and Fermat's last theorem,
4th edition, Boca Raton, FL: CRC Press, 2016.

\end{thebibliography}
\end{document}